\RequirePackage{fix-cm}
\documentclass[smallextended]{svjour3}       
\journalname{Advances in Statistical Analysis}
\smartqed  
\usepackage{graphicx}
\usepackage{amssymb}
\usepackage{natbib}
\usepackage{amsmath,enumerate,bbm,color,verbatim,setspace}
\begin{document}

\title{Nonparametric estimation of Multivariate density and its
derivative  by Dependent data using gamma kernels
}

\titlerunning{Multivariate gamma kernel estimation by dependent data}        

\author{Liubov A. Markovich        
}


\institute{Liubov A. Markovich $^{1,2,3}$ \at
      $^1$Institute of Control Sciences, RAS, Moscow\\
$^2$Moscow Institute of Physics and Technology,\\
$^3$Institute for information transmission problems, RAS, Moscow\\
              \email{kimo1@mail.ru }           
}

\date{Received: date / Revised: date}

\maketitle

\begin{abstract}
We consider the nonparametric estimation
 of the multivariate probability density function and its partial derivative with a support on
 $[0,\infty)$ by dependent data. To this end we use the class of kernel estimators with asymmetric
gamma kernel functions. The gamma kernels are nonnegative. They change their shape depending on the position
on the semi-axis and possess good boundary properties for a wide class of densities.
The theoretical asymptotic properties  of the multivariate density and its partial derivative estimates like biases, variances and covariances are derived.
We obtain the optimal bandwidth selection for both estimates as a minimum of the mean integrated squared error (MISE)
assuming dependent data with a strong mixing. Optimal rates of convergence of the MISE both for the density and its derivative are found.
\end{abstract}

\section{Introduction}
Nonnegatively supported probability density functions (pdf) can be used to model a wide range of applications in engineering, signal processing, medical research, quality control, actuarial science and climatology among others. Regarding the optimal filtering in the signal processing and control of nonlinear processes the exponential pdf class is used (e.g, \citet{Dobrovidov:12}). Most total insurance claim distributions are shaped by gamma pdfs (cf. \citet{Furman}) that are nonnegatively supported, skewed to the right and unimodal. The gamma distributions are also used to model rainfalls (e.g, \citet{Aksoy}). Erlang and $\chi^2$ pdfs are widely used in modeling insurance portfolios (e.g, \citet{Hurlimann}).
\par The wide use of these pdfs in practice leads to the need of their estimation by data samples of moderate sizes. Kernel estimators seem to be one of the most popular nonparametric  estimators.
The kernel density estimators were originally introduced for univariate independent identically distributed (iid) data and
for symmetrical kernels in \citet{Parzen,Rosenblatt}. The latter approach was widely adopted for tasks where the support of the underlying pdf $f$ is unbounded. In case when the pdf has a nonsymmetric support $[0,\infty)$ the problem of a large bias on the zero boundary appears. This leads to a bad quality of the estimates (cf. \citet{WandJones}). Obviously, the boundary bias  for the multivariate pdf estimation becomes even more solid (e.g, \citet{TaufikBouezmarnia:Rom}).
To overcome this problem one can use special approaches such as the data reflection (e.g. \citet{Schuster}), the boundary kernels (e.g, \citet{Muller}), the hybrid method (e.g. \citet{HallWehrly}), the local linear estimator (e.g, \citet{LejeuneSarda,Jones}) among others. Another solution  is to use asymmetrical kernel functions instead of  symmetrical ones.
For univariate nonnegative iid random variables (r.v.s), the  estimators with gamma kernels were proposed in \citet{Chen:20}.  The gamma kernel estimator was developed for univariate dependent data in \citet{TaufikBouezmarnia:Rom2}. In \citet{TaufikBouezmarnia:Rom} the gamma kernel estimator of the multivariate pdf for nonnegative iid r.v.s was introduced.
\par Gamma kernel is nonnegative and flexible regarding the shape. This allows to fit accurately multi-modal pdfs and their derivatives. Estimators constructed with gamma kernels have no boundary bias if $f''(0)=0$ holds, i.e when the underlying pdf $f(x)$ has a shoulder at $x = 0$ (see formula (4.3) in \citet{Zhang}). This shoulder property is fulfilled, for example, for a wide exponential class of pdfs. Several bias correction methods without the shoulder property can be found in \citet{Igarashi} for univariate iid data  and in \citet{Funke} for multivariate iid data. Other asymmetrical kernel estimators like inverse Gaussian and reciprocal inverse Gaussian estimators were studied in \citet{Scailet}. The comparison of these asymmetric kernels with the gamma kernel is given in \citet{BouSca}.
\par Along with the pdf estimation it is often necessary  to estimate the derivative of the pdf. The estimation of the univariate pdf derivative by the gamma kernel estimator was proposed  in \citet{DobrovidovMarkovich:13a,DobrovidovMarkovich:13b} for iid data and in \citet{Mar2} for a strong mixing dependent data. Our procedure achieves the optimal MISE of order  $n^{-4/7}$ when the optimal bandwidth is of order $n^{-2/7}$.
In \citet{WandJones} an optimal MISE of the kernel estimate of the first derivative of order $n^{-4/7}$  corresponding to the optimal bandwidth of order $n^{-1/7}$ for symmetrical kernels was indicated.
\subsection{Contributions of this paper}
In this paper, we introduce for the first time the gamma product kernel estimators for the  multivariate joint pdf with a nonnegative support and its partial derivative by multivariate dependent data with a strong mixing. The asymptotical behavior of the estimates and  the optimal bandwidths in the sense of minimal MISE are obtained. Note that the derivative estimation requires a specific bandwidth selection that is different from that for the pdf estimation.
\subsection{Practical motivation}
In practice it is necessary  to deal with sequences of dependent observations. As a natural relaxation of the iid condition one can deal with stationary processes satisfying the strong mixing condition.
As an example of such processes one can take  autoregressive processes. Along with the evaluation of the density function and its derivative
by dependent samples, the estimation of the logarithmic derivative of the density is an actual problem (see the equation for the optimal filtering in \citet{Marlit}  eq. (2.4), (2.8)). The logarithmic pdf derivative  is the ratio of the derivative of the pdf to the pdf itself. The pdf derivative  estimation have practical use in an optimal filtering in signal processing and control of nonlinear processes where only the exponential pdf class is used (cf. \citet{Dobrovidov:12,Dobr}). For more details see Sec. \ref{sec:5_1}.
\par The pdf derivative could be useful to find a slope of the pdf curve, its local extremes, saddle points, significant features in data as well as in regression analysis (cf. \citet{Brabanter}). The pdf derivative plays also  a key role in data clustering by means of mode seeking (cf. \citet{Sasaki}).
\par The outline of the paper is as follows. In Section \ref{lab:1} we define the gamma kernel estimators of the multivariate joint pdf and its partial derivative. In Section \ref{lab:2} we obtain the bias, the variance and the covariance of the joint pdf estimate. Using these results we derive the  optimal bandwidth and the corresponding rate of the optimal MISE. In Section \ref{lab:3} we obtain the same for the estimate of the partial pdf derivative. In Section \ref{lab:4} we investigate the moderate sample properties of the  gamma kernel density estimator for
nonnegative bivariate data. Section \ref{lab:5} provides the conclusion. The proofs of the theorems are presented in the Appendix.
\section{Gamma Kernel Estimation}\label{lab:1}
Let $\{\mathbf{X}_i = (X_{i1},\ldots,X_{id})^T\}_{i=1}^{n}$ be a strongly stationary sequence of $d$-dimensional variables with an unknown pdf  $f(x_1^{d})$ such that $\{X_{ij}\}$ are iid r.v.s. Let $x_1^{d}=(x_1,\ldots,x_d)^T$ be defined at the bounded support $x_1^{d}\in \mathbf{R}^{+d}$. We assume  that the sequence $\{\mathbf{X}_{s}\}$ is $\alpha-$mixing with coefficient
\begin{equation*}
    \alpha(k)  = \sup\limits_{t}\sup\limits_{\genfrac{}{}{0pt}{}{A\in \mathfrak{\sigma}\{X_s,s\leq t\}}{B\in \mathfrak{\sigma}\{X_s,s\geq t+k\}}}|P(A\cap B)-P(A)P(B)|, \quad k\geq1.\label{1}
\end{equation*}
Here $\alpha(k)\rightarrow0$ as $k\rightarrow\infty$. The $\alpha$-mixing condition, also called strong mixing, is satisfied by many stochastic processes, e.g, by autoregressive processes
(e.g, \citet{Donald:83}).
For these sequences we will use the notation $\{\mathbf{X}_s\}\in\mathcal{S}(\alpha)$.
To estimate the unknown multivariate pdf we use the product gamma kernel estimator (see \citet{TaufikBouezmarnia:Rom})
\begin{eqnarray}\label{35a}\widehat{f}(x_1^{d})&=&\frac{1}{n}\sum\limits_{i=1}^{n}\prod\limits_{j=1}^{d}K_{{\rho(x_j.b_j),b_j}}\left(X_{ij}\right),
\end{eqnarray}
where $b\equiv\{b_{j}\}_{j=1}^{d}$ is the vector of the bandwidth parameters such that $b\rightarrow 0$ as $n\rightarrow\infty$.
We propose the gamma kernel (cf. \citet{Chen:20})
\begin{equation}\label{Krho}K_{\rho(x,b),b}(t)=\frac{t^{\rho(x,b)-1}\exp(-t/b)}{b^{\rho(x,b)}\Gamma(\rho(x,b))},
\end{equation}
as a kernel function for each variable, where $\Gamma(\cdot)$ is a standard gamma function. To improve the properties of the gamma kernel estimator in \citet{Chen:20} it was proposed to use the parameter $\rho(x,b)$ defined as
\begin{eqnarray}\label{rho}
\rho(x,b)&=& \left\{
\begin{array}{ll}
\rho_1(x,b) = x/b, &   \mbox{if}\qquad x\geq 2b,
\\
\rho_2(x,b) =\left(x/(2b)\right)^2+1, & \mbox{if}\qquad x\in
[0,2b).
\end{array}
\right.
\end{eqnarray}
Since the gamma kernel is nonnegative it has no weight on the negative semiaxes in contrast to symmetrical kernel estimators. Hence, the use of the gamma kernels is
more natural for the estimation of the nonnegatively supported pdfs.
Using \eqref{rho} the pdf estimator \eqref{35a} can be written as
\begin{eqnarray}\label{7}
\widehat{f}(x_{1}^{d})&=&\left\{
\begin{array}{ll}
\frac{1}{n}\sum\limits_{i=1}^{n}\prod\limits_{j=1}^{d}K_{{\rho_1(x_j,b_j),b_j}}\left(X_{ij}\right),
& \mbox{if}\quad x\geq 2b,\\
\frac{1}{n}\sum\limits_{i=1}^{n}\prod\limits_{j=1}^{d}K_{{\rho_2(x_j,b_j),b_j}}\left(X_{ij}\right),&
\mbox{if}\quad x\in [0,2b).
\end{array}
\right.
\end{eqnarray}
In \citet{DobrovidovMarkovich:13a,Mar2}
the derivative of the univariate pdf was estimated just like the derivative of the gamma kernel estimator. For the multivariate case, we can analogically estimate any
partial derivative of $f(x_{1}^{d})$. For example, the partial derivative of $f(x_{1}^{d})$ by $x_k$, $1\leq k\leq d$ can be estimated as
\begin{eqnarray}\label{f'(x)}
\widehat{f}'_{x_k}(x_{1}^{d})&=&\frac{1}{n}\sum\limits_{i=1}^{n}\prod\limits_{j=1}^{d}K_{{\rho(x_j,b_j),b_j}}\left(X_{ij}\right)(K_{{\rho(x_k,b_k),b_k}}\left(X_{ik}\right))'_{x_k},
\end{eqnarray}
where
\begin{eqnarray}\label{K'}
(K_{\rho(x,b),b}(t))'_{x}&=&\left\{
\begin{array}{ll}
K'_{\rho_1(x,b),b}(t)=\frac{1}{b}K_{\rho_1(x,b),b}(t)L_1(t,x,b),
& \mbox{if}\quad \!\!x\geq 2b,\\
K'_{\rho_2(x,b),b}(t)=\frac{x}{2b^2}K_{\rho_2(x,b),b}(t)L_2(t,x,b),&
\mbox{if}\quad\!\! x\in [0,2b),
\end{array}
\right.
\end{eqnarray}
is the partial derivative of \eqref{Krho} and
\begin{eqnarray*}\label{L}
L_i(t,x,b)&=& \ln t - \ln b - \Psi(\rho_i(x,b)),\quad i=1,2.
\end{eqnarray*}
Here, $\Psi(\rho_i(x,b))$  denotes  the Digamma function that is the logarithmic derivative of the gamma
function.
\begin{remark}The mathematical tool applied for the derivative estimation is similar to one applied for the pdf. However, formulas become much more complicated
because one has to deal with the special Digamma function that includes the bandwidth vector $b$. Hence, one has to pick out the order by $b$ from complicated expressions containing logarithms and a special function.
\end{remark}
Using \eqref{K'}, formula \eqref{f'(x)} can be rewritten by
\begin{eqnarray}\label{6}
\widehat{f}'_{x_k}(x_{1}^{d})&=&\left\{
\begin{array}{ll}
\frac{1}{n}\sum\limits_{i=1}^{n}\frac{1}{b_k}L_1(X_{ik},x_k,b_k)\prod\limits_{j=1}^{d}K_{{\rho_1(x_j,b_j),b_j}}\left(X_{ij}\right),
& \mbox{if}\quad x\geq 2b,\\
\frac{1}{n}\sum\limits_{i=1}^{n}\frac{x_k}{2b_k^2}L_2(X_{ik},x_k,b_k)\prod\limits_{j=1}^{d}K_{{\rho_2(x_j,b_j),b_j}}\left(X_{ij}\right),&
\mbox{if}\quad x\in [0,2b).
\end{array}
\right.
\end{eqnarray}
It is natural to consider the MISEs which are defined as
\begin{eqnarray}\label{3a}
MISE(\widehat{f}(x_{1}^{d}))&=&
\mathsf E\int\limits_0^\infty(f(x_{1}^{d})-\widehat{f}(x_{1}^{d}))^2dx, \\\nonumber
MISE(\widehat{f}'_{x_k}(x_{1}^{d}))&=&
\mathsf E\int\limits_0^\infty(f'_{x_k}(x_{1}^{d})-\widehat{f}'_{x_k}(x_{1}^{d}))^2dx
\end{eqnarray}
as the measure of error of the proposed estimators \eqref{7} and \eqref{6}. The unknown smoothing parameters of \eqref{7} and \eqref{6} are obtained as minima of \eqref{3a}.
\begin{remark}The integrals in \eqref{3a} can be splitted into two integrals $\int_{0}^{2b}$ and $\int_{2b}^{\infty}$.
Further we shall do all the proofs for the case when $x\geq 2b$ because the integral $\int_{0}^{2b}$ tends to zero if $b\rightarrow0$. Hence, we omit the indices for $\rho_1(x_j,b_j)$, $L_1(X_{ik},x_k,b_k)$ and instead use  $\rho(x_j,b_j)$, $L(X_{ik},x_k,b_k)$.
\end{remark}
\par The mean squared error (MSE)  is determined as
\begin{eqnarray}\label{MSE}MSE(\widehat{f}(x_{1}^{d}))&=& (Bias(\widehat{f}(x_{1}^{d}))^2+Var(\widehat{f}(x_{1}^{d})),\\\nonumber
MSE(\widehat{f}'_{x_k}(x_{1}^{d}))&=& (Bias(\widehat{f}'_{x_k}(x_{1}^{d}))^2+Var(\widehat{f}'_{x_k}(x_{1}^{d})),
\end{eqnarray}
where  the variances contain for dependent r.v.s  the covariance terms
\begin{eqnarray}\label{22}&&Var(\widehat{f}(x_{1}^{d}))=
Var\left(\frac{1}{n}\sum\limits_{i=1}^{n}\prod\limits_{j=1}^{d}K_{{\rho(x_j,b_j),b_j}}\left(X_{ij}\right)\right)\\\nonumber
&=&\frac{1}{n^2}\Biggl(\sum\limits_{i=1}^{n}Var\left(\widetilde{\widetilde{K}}\left(X_{i},x,b\right)\right)
+\sum\limits_{i,j=1,i\neq j}^{n}\!\!\!\!Cov\left(\widetilde{\widetilde{K}}\left(X_{i},x,b\right),\widetilde{\widetilde{K}}\left(X_{j},x,b\right)\right)\Biggr)\\\nonumber
&=&\frac{1}{n}\overline{Var}(\widehat{f}(x_{1}^{d}))
+\frac{2}{n}\sum\limits_{i=1}^{n}\left(1-\frac{i}{n}\right)Cov\left(\widetilde{\widetilde{K}}\left(X_{i},x,b\right),\widetilde{\widetilde{K}}\left(X_{1+i},x,b\right)\right),
\end{eqnarray}
\begin{eqnarray}\label{23}&&Var(\widehat{f}'_{x_k}(x_{1}^{d}))=
Var\left(\frac{1}{n}\sum\limits_{i=1}^{n}\frac{1}{b_n}L(X_{in},x_n,b_n)\prod\limits_{j=1}^{d}K_{{\rho(x_j,b_j),b_j}}\left(X_{ij}\right)\right)\\\nonumber
&=&\frac{1}{n^2b_k^2}\Biggl(\sum\limits_{i=1}^{n}Var\left(\widetilde{K}\left(X_{i},x,b\right)\right)
+\sum\limits_{i,j=1,i\neq j}^{n}\!\!\!\!Cov\left(\widetilde{K}\left(X_{i},x,b\right),\widetilde{K}\left(X_{j},x,b\right)\right)\Biggr)\\\nonumber
&=&\frac{1}{nb_k^2}\overline{Var}(\widehat{f}'_{x_k}(x_{1}^{d}))
+\frac{2}{nb_k^2}\sum\limits_{i=1}^{n}\left(1-\frac{i}{n}\right)Cov\left(\widetilde{K}\left(X_{1},x,b\right),\widetilde{K}\left(X_{1+i},x,b\right)\right).
\end{eqnarray}
Here and below we use the following notations
\begin{eqnarray*}\widetilde{\widetilde{K}}\left(X_{i},x,b\right)=\prod\limits_{j=1}^{d}K_{{\rho(x_j,b_j),b_j}}\left(X_{ij}\right),\quad\!\!
\widetilde{K}\left(X_{i},x,b\right)=L(X_{ik},x_k,b_k)\widetilde{\widetilde{K}}\left(X_{i},x,b\right).
\end{eqnarray*}
\section{Main results}\label{lab:2}
\subsection{Convergence rate of the density estimator}
In this section we obtain the asymptotic properties of the estimator \eqref{7}. To this end we derive the bias, the variance and the covariance determined in \eqref{MSE} in the following lemmas. All proofs are relegated to the Appendix and hold assuming that all the components of the bandwidth vector are different. Regarding the practical use we give a simpler formulation of the lemmas for equal bandwidths $b_1=b_2=\ldots=b_d=b$.
\\In the next lemmas devoted to the bias and the variance we assume that $f(x_{1}^d)$ is a twice continuously differentiable function. The next lemma states the bias of the nonparametric density estimator.
\begin{lemma}\label{lem3}If $b_1=b_2=\ldots=b_d=b$  and $b\rightarrow 0$ hold as $n\rightarrow \infty$, then the
bias of the pdf estimate  \eqref{7} is equal to
\begin{eqnarray}\label{28}
Bias(\widehat{f}(x_{1}^{d}))&=&
\frac{b}{2}\sum\limits_{j=1}^{d}x_j\frac{\partial^2 f(x_{1}^{d})}{\partial x_j^2}+o\left(b\right).
\end{eqnarray}
\end{lemma}
The following result states the variance of the nonparametric estimator.
\begin{lemma}\label{lem4}If $b_1=b_2=\ldots=b_d=b$ and $b\rightarrow 0$, $
nb^{\frac{d}{2}}\rightarrow \infty$ hold as $n\rightarrow \infty$, then the
 variance expansion of the pdf estimate
\eqref{7} is equal to
\begin{eqnarray}\label{29}&&Var(\widehat{f}(x_{1}^{d}))
=\frac{b^{-\frac{d}{2}}}{n}\left(\prod\limits_{j=1}^{d}\frac{x_j^{-1/2}}{2\sqrt{\pi}}\right)\left(f(x_{1}^d)+bv_1(x_{1}^d)+b^2v_2(x_{1}^d)\right)\nonumber\\
&-&\frac{1}{n}\left(f(x_{1}^{d})+\frac{1}{2}\sum\limits_{j=1}^{d}x_jb\frac{\partial^2 f(x_{1}^{d})}{\partial x_j^2}\right)^2+o(b^2),
\end{eqnarray}
where
\begin{eqnarray*}\label{30}v_1(x_{1}^d)&=&\sum\limits_{j=1}^{d}\left(-\frac{1}{2}\frac{\partial f(x_{1}^d)}{\partial x_j}+\frac{x_j}{4}\frac{\partial^2 f(x_{1}^d)}{\partial x_j^2}\right),\quad
v_2(x_{1}^d)=-\sum\limits_{j=1}^{d}\sum\limits_{i=1}^{d}\frac{x_j}{8}\frac{\partial^3 f(x_{1}^d)}{\partial x_j^2\partial x_i}.
\end{eqnarray*}
\end{lemma}
Now we turn our attention to the covariance introduced in \eqref{22}.
\begin{lemma}\label{lem7}Let
\begin{enumerate}
  \item $\{X_j\}_{j\geq1}\in\mathcal{S}(\alpha)$ and $\int\limits_1^\infty \alpha(\tau)^\upsilon d\tau<\infty,\quad 0<\upsilon<1$ hold,
  \item $f(x_{1}^d)$ is a twice continuously differentiable function,
  \item $b\rightarrow 0$ and $nb^{d\frac{\upsilon+1}{2}}\rightarrow \infty$ as $n\rightarrow\infty$.
 \end{enumerate}
  Then, the covariance is bounded by
 \begin{eqnarray*}\label{31}|Cov(\widehat{f}(x_{1}^{d})|\leq\frac{b^{-d\frac{\upsilon+1}{2}}}{n}\int\limits_{1}^{\infty}\alpha(\tau)^{\upsilon}d\tau D(\upsilon,x_{1}^d)\Bigg(bS(\upsilon,x_{1}^d)+f(x_{1}^d)\frac{3\upsilon-1}{2(\upsilon-1)}\Bigg)^{1-\upsilon}\!\!\!\!\!\!+o\left(b^2\right),
\end{eqnarray*}
where
\begin{eqnarray*}S(\upsilon,x_{1}^d)&=&\sum\limits_{i=1}^d
\frac{\upsilon+1}{(\upsilon-1)^2x_i}f(x_{1}^d)+\frac{\upsilon+1}{\upsilon-1}\frac{\partial f(x_{1}^d)}{\partial x_i}+\frac{x_i}{2}\frac{\partial^2 f(x_{1}^d)}{\partial x_i^2},
\end{eqnarray*}
\begin{eqnarray*}D(\upsilon,x_{1}^d)&=&2(2\pi)^{-\frac{d(\upsilon+1)-2}{2}}\Bigg(\prod\limits_{j=1}^{d}x_j^{-\frac{\upsilon+1}{2}}\Bigg).
\end{eqnarray*}
\end{lemma}
Using the results of the latter we can obtain the upper bound of the MISE \eqref{3a}
and find the expression of the optimal bandwidth $b$ as the minimum of the latter.
Hence, the following theorem holds.
 \begin{theorem}\label{thm2}In conditions of Lemmas \ref{lem3} - \ref{lem7},
the optimal bandwidth that provides a minimum of the MISE is given by
\begin{eqnarray}\label{11}&&b=\left(\frac{d(\upsilon+1)(3\upsilon-1)^{1-\upsilon}}{n(2\upsilon-2)^{1-\upsilon}}
\frac{\int D(\upsilon,x_{1}^d)f(x_{1}^d)^{1-\upsilon}dx_{1}^{d}}{\int \left(\sum\limits_{j=1}^dx_j\frac{\partial f(x_{1}^d)}{\partial x_j^2}\right)^2dx_{1}^{d}}\int\limits_{1}^{\infty}\alpha(\tau)^{\upsilon}d\tau\right)^{\frac{d(\upsilon+1)+4}{2}},
\end{eqnarray}
where we denote multiple integrals as $\int\ldots\int dx_{1}\ldots dx_{d}$ as $ \int dx_{1}^{d}$.
 \end{theorem}
\par Obviously, it is impossible to use the bandwidth \eqref{11} since it contains the unknown mixing coefficient $\alpha(\tau)$. The following lemma serves to overcome this problem.
\begin{lemma}\label{lem77} Let the conditions of Lemma \ref{lem7} be satisfied. If additionally
 there exists $M>0$ such that for all $\tau >1$ and $(x,y)\in\mathbb{R}^2$ it follows $$|f_\tau(x,y)-f(x)f(y)|\leq M,$$
  then the covariance of the multivariate pdf estimator is given by
 \begin{eqnarray*}&&|Cov(\widehat{f}(x_{1}^{d})|\leq\frac{2}{n}\sum\limits_{\tau=2}^{c(n)}Cov\left(\widetilde{\widetilde{K}}\left(X_{1},x,b\right),\widetilde{\widetilde{K}}\left(X_{\tau},x,b\right)\right)\\
 &+&\frac{2}{n}\sum\limits_{\tau=c(n)+1}^{\infty}Cov\left(\widetilde{\widetilde{K}}\left(X_{1},x,b\right),\widetilde{\widetilde{K}}\left(X_{\tau},x,b\right)\right)=I_1+I_2
\end{eqnarray*}
and its rate of convergence is
 \begin{eqnarray}\label{32}I_1&\leq& \frac{2M}{nb^{\frac{d}{8}}},\\\nonumber
I_2&\leq&\frac{D(\kappa,x_{1}^d)}{nb^{\frac{d}{16}}}
\Bigg(f(x_{1}^d)\frac{6\kappa-1}{2(2\kappa-1)}+bS(\kappa,x_{1}^d)\Bigg)^{1-2\kappa}\sum\limits_{\tau=2}^{\infty}\tau\alpha(\tau)^{2\kappa}+o\left(b^2\right).
\end{eqnarray}
\end{lemma}
\par Due to Lemma \ref{lem77} the following rate of the covariance $Cov(\widehat{f}_n(x_{1}^d))\sim n^{-1}b^{-\frac{d}{8}}$ holds. Thus, the covariance can be neglected in comparison  with the variance $var(\widehat{f}(x_{1}^d))\sim n^{-1}b^{-\frac{d}{2}}$. The next theorem gives the optimal bandwidth that minimizes MISE of the estimator.
\begin{theorem}\label{thm6} In conditions of Lemmas \ref{lem3}, \ref{lem4}, \ref{lem77},
the optimal bandwidth which provides a minimum of the MISE is given by
\begin{eqnarray}\label{33}&&b=\left(\frac{d\int\limits_{0}^{\infty}\left(\prod\limits_{j=1}^{d}\frac{x_j^{-1/2}}{2\sqrt{\pi}}\right)f(x_{1}^d)dx_{1}^d}{\int\limits_{0}^{\infty} \left(\sum\limits_{j=1}^{d}x_j\frac{\partial^2 f(x_{1}^{d})}{\partial x_j^2}\right)^2dx_{1}^d}\right)^{\frac{2}{4+d}}n^{-\frac{2}{4+d}}.
\end{eqnarray}
 \end{theorem}
\begin{corollary}The result of Theorem \ref{thm6} is in agreement with the results of \citet{Chen:20} for univariate pdf and the iid case.
\end{corollary}
\subsection{Convergence rate of the density derivative estimator}\label{lab:3}
Similarly to Section \ref{lab:2} we obtain the asymptotic properties of the pdf derivative estimator \eqref{6}. The bias and the variance determined in \eqref{MSE} are derived  in the following lemmas.
\begin{lemma}\label{lem1}If $b_1=b_2=\ldots=b_d=b$ and  $b\rightarrow 0$ as $n\rightarrow \infty$ hold, then the
bias for the pdf derivative
estimate \eqref{6} is equal to
\begin{eqnarray}\label{34}
Bias(\widehat{f}'_{x_k}(x_{1}^{d}))
&=&bB_1(x_{1}^d)+b^2B_2(x_{1}^d)+o\left(b\right),
\end{eqnarray}
where we denote
\begin{eqnarray*}B_1(x_{1}^d)&=&\frac{f(x_{1}^{d})}{12x_k^2}+\frac{1}{4x_k}\sum\limits_{j=1}^{d}x_j\frac{\partial^2 f(x_{1}^{d})}{\partial x_j^2},\quad
B_2(x_{1}^d)=\frac{1}{24x_k^2}\sum\limits_{j=1}^{d}x_j\frac{\partial^2 f(x_{1}^{d})}{\partial x_j^2}.
\end{eqnarray*}
\end{lemma}
\begin{lemma}\label{lem2}If $b_1=b_2=\ldots=b_d=b$ and $b\rightarrow 0$, $nb^{\frac{d}{2}}\rightarrow \infty$ as $n\rightarrow \infty$ hold, then the
 variance of the pdf derivative estimate
\eqref{6} is equal to
\begin{eqnarray}\label{38}&&\!\!\!\!\!\!Var(\widehat{f}'_{x_k}(x_{1}^{d}))=\frac{b^{-\frac{d}{2}}}{n}\left(\prod\limits_{j=1}^{d}\frac{x_j^{-1/2}}{2\sqrt{\pi}}\right)
\Bigg(bV_1(x_{1}^d)+b^2V_2(x_{1}^d)
+\frac{1}{b}V_3(x_{1}^d)+V_4(x_{1}^d)\Bigg)\nonumber\\
&-&\frac{1}{n}\Bigg(b^2B_1^2(x_{1}^d)+\left(\frac{\partial f(x_{1}^{d})}{\partial x_k}\right)^2
+2\frac{\partial f(x_{1}^{d})}{\partial x_k}\left(bB_1(x_{1}^d)+b^2B_2(x_{1}^d)\right)\Bigg)\!\!\!\!+o\left(b^2\right),
\end{eqnarray}
where
\begin{eqnarray*}V_1(x_{1}^d)&=&-\frac{1}{24x_k^2}\frac{\partial f(x_{1}^d)}{\partial x_k}+\sum\limits_{j=1}^{d}\left(\frac{1}{8x_k}\frac{\partial^2 f(x_{1}^d)}{\partial x_j^2}-\frac{1}{8x_k^2}\frac{\partial f(x_{1}^d)}{\partial x_j}\right)
+\frac{7}{48x_k^3}f(x_{1}^d),
\end{eqnarray*}
\begin{eqnarray*}V_2(x_{1}^d)&=&\frac{7}{576x_k^4}f(x_{1}^d)+\sum\limits_{j=1}^{d}\left(\frac{1}{16x_k^2}\frac{\partial^2 f(x_{1}^d)}{\partial x_j^2}
-\frac{7}{96x_k^3}\frac{\partial f(x_{1}^d)}{\partial x_j}+\frac{1}{48x_k^2}\frac{\partial^2 f(x_{1}^d)}{\partial x_k\partial x_j}\right),
\end{eqnarray*}
\begin{eqnarray*}\label{37}V_3(x_{1}^d)&=&\frac{f(x_{1}^d)}{2x_k},\quad V_4(x)=\frac{f(x_{1}^d)}{4x_k^2}-\sum\limits_{j=1}^{d}\frac{1}{4x_k}\frac{\partial f(x_{1}^d)}{\partial x_j}.
\end{eqnarray*}
\end{lemma}
Finally,we find the covariance determined in \eqref{23} and its rate of convergence in the next lemma .
\begin{lemma}\label{lem5}Let
\begin{enumerate}
  \item $\{X_j\}_{j\geq1}\in\mathcal{S}(\alpha)$ and $\int\limits_1^\infty \alpha(\tau)^\upsilon d\tau<\infty,\quad 0<\upsilon<1$ hold,
  \item $b\rightarrow 0$ and $nb^{d\frac{\upsilon+1}{2}}\rightarrow \infty$ as $n\rightarrow\infty$.
 \end{enumerate}
  Then the covariance of the estimate of the pdf partial derivative is bounded by
   \begin{eqnarray*}\label{cov}|Cov(\widehat{f}'_{x_k}(x_{1}^{d})|&\leq&
   \frac{R(\upsilon,x_{1}^d)}{nb^{d\frac{\upsilon+1}{2}}}\Bigg(b^2V(\upsilon,x_{1}^d)+bW(\upsilon,x_{1}^d)
   +L(\upsilon,x_{1}^d)\Bigg)^{1-\upsilon}\nonumber\\
   &\cdot&\int\limits_{1}^{\infty}\alpha(\tau)^{\upsilon}d\tau+o(b^2),
\end{eqnarray*}
where we denote
 \begin{eqnarray*}V(\upsilon,x_{1}^d)&=&\sum\limits_{i=1}^d\bigg(\bigg(\frac{(\upsilon+1)(3\upsilon-1)}{72(\upsilon-1)^3x_k^2}
 +\frac{\upsilon+1}{(\upsilon-1)^2x_i}-\frac{\upsilon(\upsilon+1)}{9(\upsilon-1)^4x_k^3}\bigg)f(x_{1}^d)\\
 &+&\frac{\upsilon+1}{\upsilon-1}\frac{\partial f(x_{1}^d)}{\partial x_i}-\frac{\upsilon(\upsilon+1)}{9(\upsilon-1)^3x_k^2}\frac{\partial f(x_{1}^d)}{\partial x_k}+
\frac{x_i}{2}\frac{\partial^2 f(x_{1}^d)}{\partial x_i^2}
\Bigg),
 \end{eqnarray*}
 \begin{eqnarray*}W(\upsilon,x_{1}^d)&=&\sum\limits_{i=1}^d\bigg(\bigg(\frac{3\upsilon-1}{4(\upsilon-1)x_k}+\frac{\upsilon+1}{(\upsilon-1)^2x_i}+\frac{2x_i(\upsilon+1)}{3(\upsilon-1)^3x_k^2}\bigg)f(x_{1}^d)\\
&+&\frac{\upsilon+1}{\upsilon-1}\frac{\partial f(x_{1}^d)}{\partial x_i}+\frac{2(\upsilon+1)x_i}{3(\upsilon-1)^2x_k}\frac{\partial f(x_{1}^d)}{\partial x_k}+\frac{x_i}{2}\frac{\partial^2 f(x_{1}^d)}{\partial x_i^2}\bigg),
 \end{eqnarray*}
  \begin{eqnarray*}L(\upsilon,x_{1}^d)=\sum\limits_{i=1}^d\Bigg(f(x_{1}^d)\frac{3\upsilon-1}{2(\upsilon-1)}+x_i\bigg(-\frac{4}{\upsilon-1}\frac{\partial f(x_{1}^d)}{\partial x_k}-\frac{4}{(\upsilon-1)^2x_k}f(x_{1}^d)\bigg)\Bigg),
   \end{eqnarray*}
  \begin{eqnarray*}R(\upsilon,x_{1}^d)&=&\left(\prod\limits_{j=1}^{d}x_j^{-\frac{\upsilon+1}{2}}\right)\frac{(2\pi)^{-\frac{d(\upsilon-3)+2}{2}}}{2x_k^2}.
   \end{eqnarray*}
\end{lemma}
\par Using the upper bound of the covariance $|C(\widehat{f}'_{x_k}(x_{1}^{d})|$ we can obtain the upper bound of the MISE and find the optimal bandwidth $b$.
The following theorem establishes the optimal bandwidth and the MISE's rate of convergence.
 \begin{theorem}\label{thm}Given the conditions of Lemmas \ref{lem1} - \ref{lem5},
the optimal bandwidth of the estimator of the partial pdf derivative by $x_k$, $1\leq k \leq d$ that provides a minimum of the MISE is given by
\begin{eqnarray}\label{35}b^{*}_{k}&=&\left(\frac{d+2}{2^{d-1}\pi^{\frac{d}{2}}}\frac{\int\limits_{0}^{\infty}\frac{f(x_{1}^d)}{x_k}\prod\limits_{j=1}^{d}x_j^{-1/2}dx_{1}^{d}}
{\int\limits_{0}^{\infty}\left(\frac{f(x_{1}^d)}{3x_k^2}+\frac{1}{x_k}\sum\limits_{i=1}^{d}x_i\frac{\partial^2 f(x_{1}^d)}{\partial x_i^2}\right)^2dx_{1}^{d}}\right)^{\frac{2}{d+6}}n^{-\frac{2}{d+6}}
\end{eqnarray}
and the corresponding MISE has the order $MISE\sim n^{-\frac{4}{d+6}}$.
 \end{theorem}
\begin{corollary}The result of Theorem \ref{thm} coincides with the results of \citet{Mar2} for the univariate pdf derivative.
\end{corollary}
\section{Application}\label{lab:4}
\par The optimal smoothing parameters for the multivariate pdf and its derivative estimate, defined by formulas \eqref{33} and \eqref{35}, depend on the
unknown true density $f(x_{1}^d)$ and its derivatives. Therefore, it is impossible to calculate values of these parameters.
\par A lot of methods for the bandwidth estimation by the sample $X_1^n$ are known. The simplest
and most convenient  is  given by the rule of thumb method (cf. \citet{Silverman:1986}) proposed for symmetrical kernels.
Instead of the unknown density $f(x_{1}^{d})$ in the optimal formulas for the bandwidth parameters a so called
reference function $g(x_{1}^{d})$ is substituted, i.e. a density in the form of the kernel function.
\par In other words, we get the estimates of optimal bandwidths by the sample in the following form
\begin{eqnarray*}\widehat{b}&=&C_1 R_1(\widehat{\mu},\widehat{\sigma}^2)n^{-\frac{2}{d+4}},\quad \widehat{b}^{*}=C_2 R_2(\widehat{\mu},\widehat{\sigma}^2)n^{-\frac{2}{d+6}},
\end{eqnarray*}
where $R_1(\widehat{\mu},\widehat{\sigma}^2)$ and $R_2(\widehat{\mu},\widehat{\sigma}^2)$ are functions depending on consistent estimates of the sample mean and variance.
The derivation of the rate by $n$ of the latter functions is out of scope of the paper.
\par For the kernel estimation with symmetrical kernels the Gaussian density is usually used as the reference function. In our paper, we use asymmetrical gamma densities.
\subsection{Bivariate gamma density}
Bivariate densities are often required for the density estimation in practice. As an example, let us construct the rule of thumb method with a bivariate reference function.
There exist many bivariate gamma-like density functions. Some of them can be found in \citet{Nadarajah1,Nadarajah2}. We select the McKay bivariate gamma distribution
whose pdf is defined as
\begin{eqnarray}\label{32_1}f(x,y)=\frac{x^{\alpha-1}(y-x)^{\beta-1}\exp(-y/\mu)}{\mu^{\alpha+\beta}\Gamma(\alpha)\Gamma(\beta)}, \quad x\leq y,\quad y>0,
\end{eqnarray}
where $\alpha,\beta>0$, $0<\mu<\infty$ and $X\sim\Gamma(\alpha,\mu)$, $Y\sim\Gamma(\alpha+\beta,\mu)$. Let us take \eqref{32_1} as the reference function $g(x,y)$. Using the method of moments we can write
\begin{eqnarray*}\overline{m}_X=\alpha\mu,\quad D_X=\alpha\mu^2,\quad\overline{m}_Y=(\alpha+\beta)\mu,\quad D_Y=(\alpha+\beta)\mu^2,
\end{eqnarray*}
where $\overline{m}_X=\frac{1}{n}\sum\limits_{i=1}^{n}X_{i}$, $\overline{m}_Y=\frac{1}{n}\sum\limits_{i=1}^{n}Y_{i}$, $\overline{D}_X=\frac{1}{n}\sum\limits_{i=1}^{n}(X_{i}-\overline{m}_X)^2$,
$\overline{D}_Y=\frac{1}{n}\sum\limits_{i=1}^{n}(Y_{i}-\overline{m}_Y)^2$.
Then, the parameters of the reference function are determined by
\begin{eqnarray*}\widehat{\mu}=\frac{\overline{D}_X}{\overline{m}_X},\quad \widehat{\alpha}=\frac{\overline{m}_X^2}{\overline{D}_X},\quad \widehat{\beta}=\frac{\overline{m}_Y\overline{m}_X}{\overline{D}_X}-\frac{\overline{m}_X^2}{\overline{D}_X}.
\end{eqnarray*}
Since the bandwidth parameter \eqref{33} contains also the second derivative of the pdf we have to use
\begin{eqnarray*}\frac{\partial^2 g(x,y)}{\partial x^2}&=&\frac{x^{\alpha-3}e^{-\frac{y}{\mu}}(y-x)^{\beta-3}}{\Gamma(\alpha)\Gamma(\beta)\mu^{\alpha+\beta}}(x^2(\alpha+\beta-2)(\alpha+\beta-3)\\
&+&y^2(\alpha-1)(\alpha-2)-2yx(\alpha-1)(\alpha+\beta-3)),\\
\frac{\partial^2 g(x,y)}{\partial y^2}&=&\frac{x^{\alpha-1}e^{-\frac{y}{\mu}}(y-x)^{\beta-3}}{\Gamma(\alpha)\Gamma(\beta)\mu^{\alpha+\beta+2}}(x^2+y^2+2x(\mu \beta-y-\mu)\\
&+&2y\mu(1-\beta)+\mu^2(\beta^2-3\beta+2)).
\end{eqnarray*}
Hence, we can take the following estimate
\begin{eqnarray}\label{34_1}&&\widehat{b}=\left(\frac{2\int\limits_{0}^{\infty}\int\limits_{0}^{\infty}\left(\frac{(x_1x_2)^{-1/2}}{4\pi}\right)\widehat{g}(x_1,x_2)dx_{1}dx_{2}}{\int\limits_{0}^{\infty}\int\limits_{0}^{\infty} \left(x_1\frac{\partial^2 \widehat{g}(x_1,x_2)}{\partial x_1^2}+x_2\frac{\partial^2 \widehat{g}(x_1,x_2)}{\partial x_2^2}\right)^2dx_1dx_2}\right)^{\frac{2}{7}}n^{-\frac{2}{7}}
\end{eqnarray}
as the optimal bandwidth parameter \eqref{33} and substitute it into the bivariate gamma kernel estimator
\begin{eqnarray}\label{33_1}\widehat{f}(X_1,X_2)&=&\frac{1}{n}\sum\limits_{i=1}^{n} \frac{X_{1i}^{x_1/\widehat{b}-1}\exp(-X_{1i}/\widehat{b})}{\widehat{b}^{x_1/\widehat{b}}\Gamma(x_1/\widehat{b})}
\frac{X_{2i}^{x_2/\widehat{b}-1}\exp(-X_{2i}/\widehat{b})}{\widehat{b}^{x_2/\widehat{b}}\Gamma(y/\widehat{b})}.
\end{eqnarray}
\subsection{Multivariate gamma density}
Let $Y=(Y_0,Y_1,\ldots,Y_n)^T$ be a $n+1$-dimensional vector of iid r.v.s, $Y_i\sim \Gamma\left(\frac{x_i}{b_i},\frac{1}{b_i}\right)$, $i=0,1,\ldots,n$.
Let us define the following matrix
\begin{eqnarray*}A = \left(
                       \begin{array}{cccccc}
                         \alpha_0/\alpha_1 & 1 & 0 & 0 & \ldots & 0\\
                          \alpha_0/\alpha_2 & \alpha_1/\alpha_2 & 1 & 0 & \ldots & 0 \\
                         \alpha_0/\alpha_3 & \alpha_1/\alpha_3 & \alpha_2/\alpha_3& 1 &\ldots & 0 \\
                         \vdots & \vdots & \vdots & \vdots & \ddots& \vdots\\
                         \alpha_0/\alpha_n & \alpha_1/\alpha_n & \alpha_2/\alpha_n & \alpha_3/\alpha_n & \ldots & 1 \\
                       \end{array}
                     \right),
\end{eqnarray*}
where $\alpha=(\alpha_1,\alpha_2,\ldots,\alpha_n)^T=\left(\frac{1}{b_1},\frac{1}{b_2},\ldots,\frac{1}{b_n}\right)$. In practice we take all $b_i\equiv b$.
Then, the joint distribution of the r.v.s $X=AY$ is called the multivariate ladder-type gamma distribution (cf. \citet{Furman}). One can take
\begin{eqnarray*}
g(x_{1}^{d})&=&\frac{x_{1}^{\frac{x_{1}}{b}}e^{-\frac{x_d}{b}}}{\Gamma(\frac{x_{1}}{b}+1)}\prod_{i=2}^{d}\frac{(x_i-x_{i-1})^{\frac{x_i}{b}-1}}{b^{\overline{\gamma}_i}\Gamma(\frac{x_i}{b})}
\end{eqnarray*}
 as a reference function alternative to a multivariate gamma density. Here $\mathbf{\overline{\gamma}}=\left(1+\frac{x_{1}}{b}, 1+\frac{x_{1}}{b}+\frac{x_{2}}{b},\ldots,1+\frac{x_{1}}{b}+\frac{x_{2}}{b}+\cdots+\frac{x_d}{b}\right)$ is the vector of the shape parameters. As in the bivariate case the unknown parameters can be set by the method of moments.
\par The proposed  rule of thumb method allows us to estimate optimal bandwidth parameters for the multivariate pdf  and its derivative by dependent data.
One can choose any other data-driven bandwidth choice instead of the proposed one.
\subsubsection{Simulation results}
In this section we study the finite sample properties for the nonparametric gamma kernel density estimator for
bivariate data with non-negative supports. To this end we use the multivariate gamma density introduced in \cite{Mathal}.
Suppose $V_1,\ldots, V_k$ are mutually independent, where $V_i\sim\Gamma(\alpha_i,\beta,\gamma_i)$, $i=1,\ldots,k$. Let 
\begin{eqnarray}\label{41}
Z_1=V_1,\quad Z_2=V_1+V_2,\ldots,\quad Z_k=V_1+\cdots +V_k.
\end{eqnarray}
Then, the joint distribution of $Z=(Z_1,\ldots,Z_k)$ is a multivariate gamma pdf (see Theorem 1.1. in \cite{Mathal}). 
One can rewrite \eqref{41} as following
\begin{eqnarray*}
Z_t=Z_{t-1}+V_t,\quad t=1,\ldots,k
\end{eqnarray*}
that is the $AR(1)$ process. According to \cite{Donald:83} it satisfies the $\alpha$-mixing condition and hence, we have an inter-temporal dependence structure.
\par We generate two i.i.d samples $V_i\sim\Gamma(\alpha_i,\mu,0)$, $i=1,2$ with sample sizes $n\in\{100, 500, 1000, 2000 \}$ using standard Matlab generators. Thus, two rv's
$Z_1=V_1$ and $Z_2=Z_1+V_2$ are dependent and their joint pdf is the McKay  bivariate gamma pdf.
\par
\begin{figure}
    \centering
        \includegraphics[width=\textwidth]{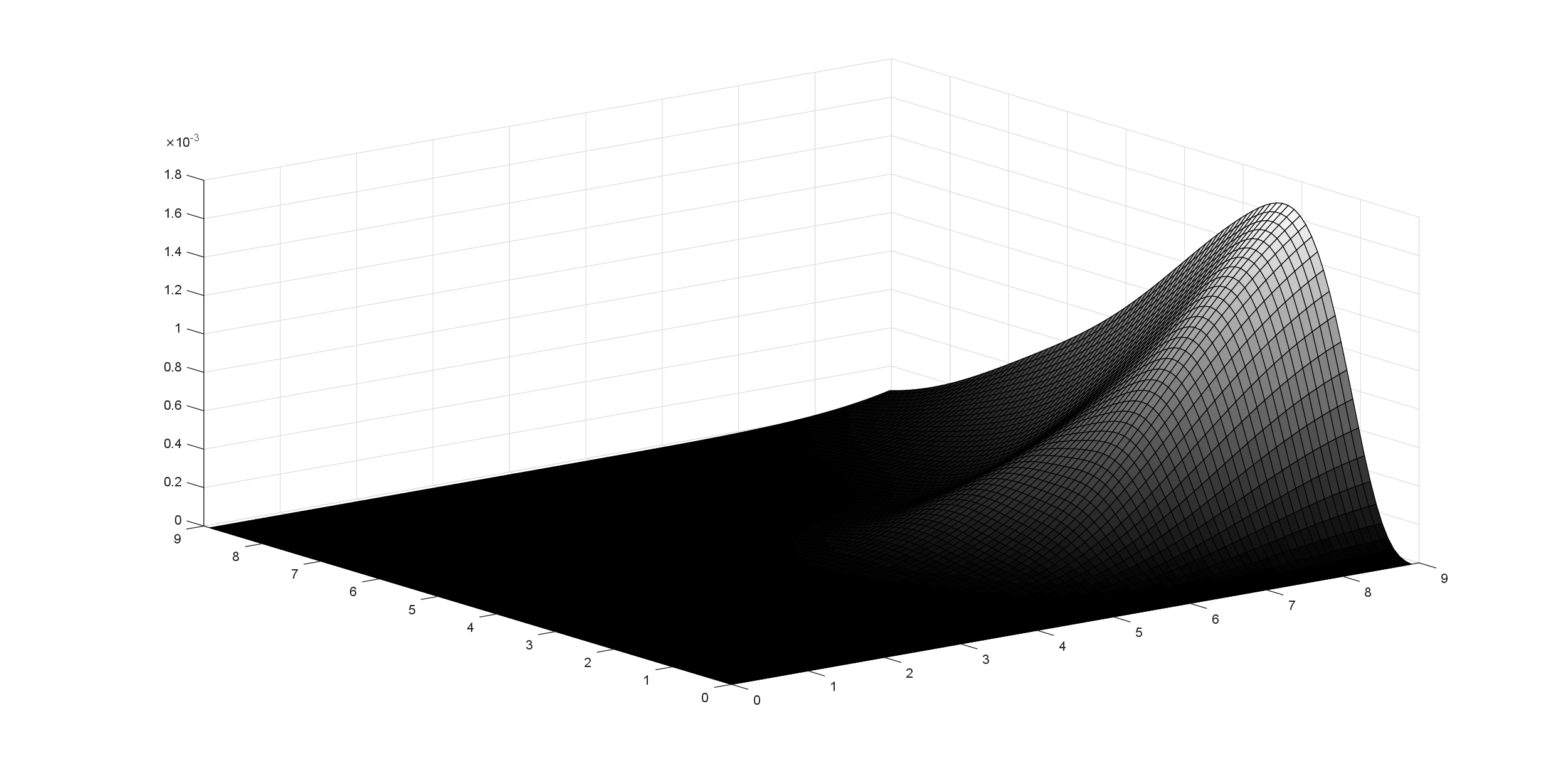}
        \caption{Estimate of the McKay pdf  for the sample size  $n=1000$}
        \label{fig:1_1}
\end{figure}

 The gamma-kernel estimate \eqref{33_1} with the optimal bandwidth is shown in Figure~\ref{fig:1_1} for $ \alpha = 3$, $\beta=5$, $\mu= 3$. The optimal bandwidth is counted for every replication of the simulation using the rule of thumb method, where  we take \eqref{32_1} as a reference density.
The estimation error of the pdf is calculated by the following formula
\begin{eqnarray*}
m=\int\limits_0^\infty(f(x,y)-\hat{f}(x,y))^2dxdy,
\end{eqnarray*}
where $f(x,y)$ is a true pdf and $\hat{f}(x,y)$ is its estimate. Values of $m$'s averaged over $500$ simulated samples and standard deviations for the underlying distribution are given in Table~\ref{Tab1}.
\begin{figure}
    \centering
        \includegraphics[width=\textwidth]{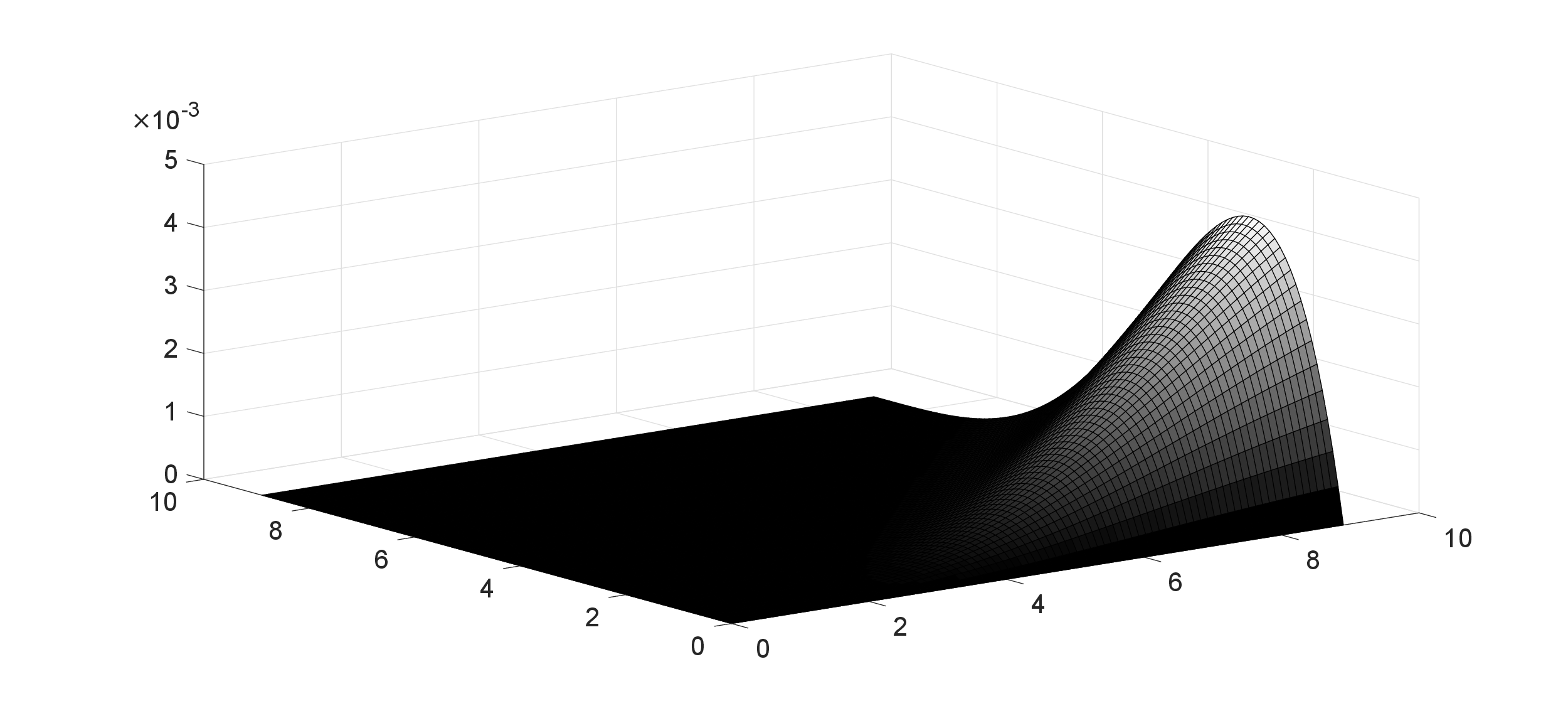}
        \caption{The  McKay pdf}
        \label{fig:2_1}
\end{figure}
\begin{table}[h]
\centering
\begin{tabular}{|c|c|c|c|c|}
  \hline
  n & 100 & 500 & 1000 & 2000 \\
  \hline
  McKay  & $1.7925 10^{-4}$  & $7.3154 10^{-5}$  &  $6.3605 10^{-5}$  &  $2.7118 10^{-5}$\\
   & $(8.1002 10^{-5})$ & $(3.8101 10^{-5})$ & $(1.8349 10^{-5})$ & $(1.2023 10^{-5})$ \\
  \hline
\end{tabular}
\caption{Mean errors and standard deviations in brackets. \label{Tab1}}
\end{table}
As expected, the mean error and the standard deviation decrease when the sample size rises.
\subsection{Practical use of the gamma kernel estimator in  optimal filtering}\label{sec:5_1}
As it was mentioned in the introduction, the results obtained in the paper have practical use in the signal processing and control of linear and nonlinear systems. The processing of stationary random sequences under nonparametric uncertainty is determined by a filtering problem when a signal distribution is unknown. A useful signal $(S_n)_{n\ge1}$ is assumed to be Markovian. This assumption allows us to estimate the unknown $(S_n)$ using only an observable random sequence $(X_n)_{n\ge1}$.
In \citet{Marlit} the equation of the optimal filtering for the exponential family of the conditional densities is given in the following form
\begin{eqnarray}\label{31}\mathsf{E}(Q(S_n)|x_1^n)\cdot T'_{x_n}(x_n)&=&\frac{f'_{x_n}(x_n|x_1^{n-1})}{f(x_n|x_1^{n-1})}-\frac{h'_{x_n}(x_n)}{h(x_n)}=
\ln\left(\frac{f(x_n|x_1^{n-1})}{h(x_n)}\right)'_{x_n}.
\end{eqnarray}
This equation contains the logarithmic derivative of the unknown conditional density which characterizes the signal. The latter density and its derivative can be written as
\begin{eqnarray*}f(x_n|x_1^{n-1})&=&\frac{f(x_1^n)}{f(x_1^{n-1})},\quad
f'_{x_n}(x_n|x_1^{n-1})=\frac{f'_{x_n}(x_1^n)}{f(x_1^{n-1})}.
\end{eqnarray*}
Hence, their ratio is determined by
\begin{eqnarray*}\frac{f'_{x_n}(x_n|x_1^{n-1})}{f(x_n|x_1^{n-1})}&=&\frac{f'_{x_n}(x_1^n)}{f(x_1^{n})}.
\end{eqnarray*}
Therefore, using the results of the previous sections we can write
\begin{eqnarray*}\frac{\widehat{f}'_{x_n}(x_{1}^n)}{\widehat{f}(x_{1}^{n})}&=&\left\{
\begin{array}{ll}
\frac{\frac{1}{b'_n}\sum\limits_{i=1}^{n}L_1(X_{in},x_n,b'_n)\prod\limits_{j=1}^{n}K_{{\rho_1(x_j,b'_j),b'_j}}\left(X_{ij}\right)}{\sum\limits_{i=1}^{n}\prod\limits_{j=1}^{n}K_{{\rho_1(x_j,b_j),b_j}}\left(X_{ij}\right)},
& \mbox{if}\quad x_{1}^n\geq 2b,\\
\frac{\frac{x_n}{2(b'_n)^2}\sum\limits_{i=1}^{n}L_2(X_{in},x_n,b'_n)\prod\limits_{j=1}^{n}K_{{\rho_2(x_j,b'_j),b'_j}}\left(X_{ij}\right)}{\sum\limits_{i=1}^{n}\prod\limits_{j=1}^{n}K_{{\rho_2(x_j,b_j),b_j}}\left(X_{ij}\right)},&
\mbox{if}\quad x_{1}^n\in [0,2b).
\end{array}
\right.
\end{eqnarray*}
In case when the bandwidth parameters are optimal, the latter estimator is the following
\begin{eqnarray*}\frac{\widehat{f}'_{x_n}(x_{1}^n)}{\widehat{f}(x_{1}^{n})}&=&\left\{
\begin{array}{ll}
\frac{\frac{1}{b'}\sum\limits_{i=1}^{n}L_1(X_{in},x_n,b')\prod\limits_{j=1}^{n}K_{{\rho_1(x_j,b'),b'}}\left(X_{ij}\right)}{\sum\limits_{i=1}^{n}\prod\limits_{j=1}^{n}K_{{\rho_1(x_j,b),b}}\left(X_{ij}\right)},
& \mbox{if}\quad x_{1}^n\geq 2b,\\
\frac{\frac{x_n}{2b'^2}\sum\limits_{i=1}^{n}L_2(X_{in},x_n,b')\prod\limits_{j=1}^{n}K_{{\rho_2(x_j,b'),b'}}\left(X_{ij}\right)}{\sum\limits_{i=1}^{n}\prod\limits_{j=1}^{n}K_{{\rho_2(x_j,b),b}}\left(X_{ij}\right)},&
\mbox{if}\quad x_{1}^n\in [0,2b).
\end{array}
\right.
\end{eqnarray*}
Next, if we substitute the latter estimate into the equation of the optimal filtering \eqref{31} we can estimate the conditional mean $\mathsf{E}(Q(S_n)|x_1^n)$
knowing only the observable part of the signal.
\section{Conclusion}\label{lab:5}
Nonparametric estimators of the probability density function and its partial derivatives on the positive semi-axis by multivariate dependent data are proposed.
Our  estimators are based on product gamma kernels with a  nonnegative support. We provide the asymptotic properties of the estimators by optimal rates of convergence
of their mean integrated squared errors. We develop explicit formulas for the optimal smoothing parameters (bandwidths) both
 for the density and the partial derivatives by samples of dependent random variables.  Furthermore, regarding positive time series the proposed gamma
kernel estimators have the same performance as Gaussian kernel estimators for the symmetrical case. Since the optimal bandwidths depend on the unknown density, it is necessary to build its databased estimate. To this end the well-known thumb method  is proposed. We examine the finite sample performance by a simulation study. Further development may concern the investigation of alternative bandwidth selection methods. The results can also be extended to samples with other mixing conditions.

\appendix
\section{Proofs}
\subsection{Proof of Lemma \ref{lem3}}\label{Ap3}
To find the bias of the estimate $\widehat{f}(x_{1}^{d})$ let us write the expectation
\begin{eqnarray}\label{12}&&\mathsf{E}_{X}(\widehat{f}(x_{1}^{d}))=
\int f(t_{1}^d)\prod\limits_{j=1}^{d}K_{\rho(x_{j},b_j)}\left(t_{j}\right)dt_{1}^{d}
=\mathsf{E}_{\xi}(f(\xi_{1}^d)),
\end{eqnarray}
 where the r.v.s $\xi_j$ are iid gamma distributed with the expectation $\mu_{j}=x_j$ and the variance
$\sigma^2_{j}=x_jb_j$. The product kernel $\prod\limits_{j=1}^{d}K_{\rho(x_{j},b_j)}\left(t_{j}\right)$ is used in \eqref{12} as a pdf.
To find $\mathsf{E}_{\xi}(f(\xi_{1}^d))$ let us do the Taylor expansion of $f(\xi_{1}^d)$  in the point $\mu_{j}$
\begin{eqnarray*}f(\xi_{1}^d)&=&f(\mu_{1}^{d})+\sum\limits_{j=1}^{d}(\xi_j-\mu_j)\frac{\partial f(\xi_{1}^d)}{\partial x_j}\Bigg|_{\xi=\mu}
+\frac{1}{2}\sum\limits_{j=1}^{d}(\xi_j-\mu_j)^2\frac{\partial^2 f(\xi_{1}^d)}{\partial x_j^2}\Bigg|_{\xi=\mu}\\
&+&\sum\limits_{j\neq l}(\xi_j-\mu_j)(\xi_l-\mu_l)\frac{\partial^2 f(\xi_{1}^d)}{\partial x_j\partial x_l}\Bigg|_{\xi=\mu}+o\left(\sum\limits_{j=1}^d b_j\right).
\end{eqnarray*}
Taking the expectation from  both sides of the latter equation we can write
\begin{eqnarray}\label{13}&&\mathsf{E}_{\xi}(f(\xi_{1}^d))=f(x_{1}^{d})+\frac{1}{2}\sum\limits_{j=1}^{d}x_jb_j\frac{\partial^2 f(x_{1}^{d})}{\partial x_j^2}+o\left(\sum\limits_{j=1}^d b_j\right).
\end{eqnarray}
Hence, the bias of the multivariate pdf estimate is given by
\begin{eqnarray*}
Bias(\widehat{f}(x_{1}^{d}))&=&
\frac{1}{2}\sum\limits_{j=1}^{d}x_jb_j\frac{\partial^2 f(x_{1}^{d})}{\partial x_j^2}+o\left(\sum\limits_{j=1}^d b_j\right).
\end{eqnarray*}
When $b_1=b_2=\ldots=b_d=b$ holds, we get \eqref{28}.
\subsection{Proof of Lemma \ref{lem4}}\label{Ap4}
By definition, the variance of the estimate $\widehat{f}(x_{1}^{d})$ is given by
\begin{eqnarray}\label{14}Var(\widehat{f}(x_{1}^{d}))
=\frac{1}{n}\left(\mathsf{E}\left(\left(\prod\limits_{j=1}^{d}K_{\rho(x_{j},b_j)}\left(X_{j}\right)\right)^2\right)
-\mathsf{E}^2\left(\prod\limits_{j=1}^{d}K_{\rho(x_{j},b_j)}\left(X_{j}\right)\right)\right).
\end{eqnarray}
Using the expectation \eqref{13}, the second term in \eqref{14} can be written as
\begin{eqnarray}\label{18}&&\frac{1}{n}\mathsf{E}_{X}^2(\widehat{f}(x_{1}^{d}))=\frac{1}{n}\left(f(x_{1}^{d})+\frac{1}{2}\sum\limits_{j=1}^{d}x_jb_j\frac{\partial^2 f(x_{1}^{d})}{\partial x_j^2}+o\left(\sum\limits_{j=1}^d b_j\right)\right)^2.
\end{eqnarray}
The first term in \eqref{14} is
\begin{eqnarray}\label{16}&&\frac{1}{n}\mathsf{E}\left(\left(\prod\limits_{j=1}^{d}K_{\rho(x_{j},b_j)}\left(X_{j}\right)\right)^2\right)
=\frac{1}{n}\int \left(\prod\limits_{j=1}^{d}K^2_{\rho(x_{j},b_j)}\left(X_{j}\right)\right)f(t_{1}^d)dt_{1}^{d}\\\nonumber
&=&\frac{1}{n}\int\left(\prod\limits_{j=1}^{d}\frac{t_j^{\frac{2x_j}{b_j}-2}\exp\left(\frac{-2t_j}{b_j}\right)}{b_j^{\frac{2x_j}{b_j}}\Gamma^2\left(\frac{x_j}{b_j}\right)}\right)
f(t_{1}^d)dt_{1}^{d}=\frac{1}{n}\left(\prod\limits_{j=1}^{d}B(x_j,b_j)\right)\mathsf{E}_{\eta}\left(f(\eta_{1}^d)\right),
\end{eqnarray}
where $\{\eta_j\}$ are iid gamma distributed r.v.s  with the expectation $\mu_{j}=x_j-\frac{b_j}{2}$ and the variance
$\sigma^2_{j}=\frac{x_jb_j}{2}-\frac{b_j^2}{4}$. We denote
\begin{eqnarray}\label{15}&&B(x_j,b_j)=\frac{b_j^{-3}x_j^2\Gamma(\frac{2x_j}{b_j}-1)}{2^{\frac{2x_j}{b_j}-1}
\Gamma^2(\frac{x_j}{b_j}+1)}.
\end{eqnarray}
Let us use the notation from \citet{Brown:99}, where $R(z) =\sqrt{2\pi}\exp(-z)z^{z+1/2}/\Gamma(z+1)$ for $z\geq 0$. Hence, we can
express the gamma functions in \eqref{15} as
\begin{eqnarray*}\Gamma^2\left(\frac{x}{b}+1\right)&=&\left(\frac{\sqrt{2\pi}\exp(\frac{-x}{b})
(\frac{x}{b})^{\frac{x}{b}+1/2}}{R(\frac{x}{b})}\right)^2,
\end{eqnarray*}
\begin{eqnarray*}\Gamma\left(\frac{2x}{b}-1\right)&=&\frac{\Gamma\left(\frac{2x}{b}+1\right)}{\frac{2x}{b}\left(\frac{2x}{b}-1\right)}=
\frac{\sqrt{2\pi}\exp\left(-\frac{2x}{b}\right)\left(\frac{2x}{b}\right)^{\frac{2x}{b}+\frac{1}{2}}}{\frac{2x}{b}\left(\frac{2x}{b}-1\right)R\left(\frac{2x}{b}\right)}.
\end{eqnarray*}
Substituting the latter expressions into \eqref{15} we can rewrite it as
\begin{eqnarray*}B(x_j,b_j)&=&
\frac{b_j^{-\frac{1}{2}}x_j^{-\frac{1}{2}}R^2(\frac{x_j}{b_j})}{2\sqrt{\pi}R(\frac{2x_j}{b_j})(1-\frac{b_j}{2x_j})}.
\end{eqnarray*}
According to Lemma 3 in \citet{Brown:99}, $R(z)$ is an
increasing function which converges to 1 as $z\rightarrow\infty$ and
$R(z)<1$ for any $z>0$.
 Then it holds
\begin{eqnarray}\label{17}
B(x_j,b_j)&=& \left\{ \begin{array}{ll}
\frac{b_j^{-1/2}x_j^{-1/2}}{2\sqrt{\pi}}, &   \mbox{if}\qquad
\frac{x_j}{b_j}\rightarrow\infty,
\\
\frac{b^{-1}k^2\Gamma(2k-1)}{2^{2k-1}\Gamma^2(k+1)}, &
\mbox{if}\qquad \frac{x_j}{b_j}\rightarrow k.
\end{array}
\right.
\end{eqnarray}
The expectation in \eqref{16} can be written similarly to \eqref{13} as
\begin{eqnarray*}&&\mathsf{E}_{\eta}(f(\eta_{1}^d))=f(\mu_{1}^d)+\frac{1}{2}\sum\limits_{j=1}^{d}\left(\frac{x_jb_j}{2}-\frac{b_j^2}{4}\right)\frac{\partial^2 f(\mu_{1}^d)}{\partial x_j^2}+o\left(\sum\limits_{j=1}^d b_j^2\right).
\end{eqnarray*}
Therefore, applying a Taylor expansion of the arguments of the functions $f(\mu_{1}^d)$
\begin{eqnarray*}f(\mu_{1}^d)&=&f(x_{1}^d)-\sum\limits_{j=1}^{d}\frac{b_j}{2}\frac{\partial f(x_{1}^d)}{\partial x_j}+\sum\limits_{j=1}^{d}\frac{b_j^2}{8}\frac{\partial^2 f(x_{1}^d)}{\partial x_j^2}+o\left(\sum\limits_{j=1}^d b_j^2\right),
\end{eqnarray*}
we deduce that \eqref{16} can be written as
\begin{eqnarray*}&&\frac{1}{n}\left(\prod\limits_{j=1}^{d}B(x_j,b_j)\right)\Bigg(f(x_{1}^d)-\sum\limits_{j=1}^{d}\frac{b_j}{2}\frac{\partial f(x_{1}^d)}{\partial x_j}+\frac{b_j^2}{8}\frac{\partial^2 f(x_{1}^d)}{\partial x_j^2}\\
&+&\frac{1}{2}\sum\limits_{j=1}^{d}\left(\frac{x_jb_j}{2}-\frac{b_j^2}{4}\right)\left(\frac{\partial^2 f(x_{1}^d)}{\partial x_j^2}-\sum\limits_{i=1}^{d}\frac{b_i}{2}\frac{\partial^3 f(x_{1}^d)}{\partial x_j^2\partial x_i}\right)\Bigg)+o\left(\sum\limits_{j=1}^d b_j^2\right).
\end{eqnarray*}
Combining \eqref{18} and \eqref{17} we can finally write the variance of the pdf estimate as
\begin{eqnarray*}&&Var(\widehat{f}(x_{1}^{d}))=\frac{1}{n}\left(\prod\limits_{j=1}^{d}\frac{b_j^{-1/2}x_j^{-1/2}}{2\sqrt{\pi}}\right)
\Bigg(f(x_{1}^d)-\sum\limits_{j=1}^{d}\frac{b_j}{2}\frac{\partial f(x_{1}^d)}{\partial x_j}\\
&+&\sum\limits_{j=1}^{d}\frac{b_j^2}{8}\frac{\partial^2 f(x_{1}^d)}{\partial x_j^2}+\frac{1}{2}\left(\frac{x_jb_j}{2}-\frac{b_j^2}{4}\right)\left(\frac{\partial^2 f(x_{1}^d)}{\partial x_j^2}-\sum\limits_{i=1}^{d}\frac{b_i}{2}\frac{\partial^3 f(x_{1}^d)}{\partial x_j^2\partial x_i}\right)\Bigg)\\
&-&\frac{1}{n}\left(f(x_{1}^{d})+\frac{1}{2}\sum\limits_{j=1}^{d}x_jb_j\frac{\partial^2 f(x_{1}^d)}{\partial x_j^2}\right)^2+o\left(\sum\limits_{j=1}^d b_j^2\right).
\end{eqnarray*}
When $b_1=b_2=\ldots=b_d=b$ holds, the latter variance can be written in the form \eqref{29}.
\subsection{Proof of Lemma \ref{lem7}}\label{Ap7}
To evaluate  the covariance of $\widehat{f}(x_{1}^{d})$ we shall apply Davydov's inequality
\begin{eqnarray*}&&|Cov\left(\widetilde{\widetilde{K}}\left(X_{1},x,b\right),\widetilde{\widetilde{K}}\left(X_{1+i},x,b\right)\right)|\leq\\
 &\cdot&2\pi \alpha(i)^{1/r}\parallel \widetilde{\widetilde{K}}\left(X_{1},x,b\right) \parallel_q \parallel \widetilde{\widetilde{K}}\left(X_{1+i},x,b\right) \parallel_p,
\end{eqnarray*}
where $p^{-1}+q^{-1}+r^{-1}=1$, $1\leq p, q, r\leq\infty$ (cf. \citet{Bosq:96}), $\alpha(i)$ is the mixing coefficient \eqref{1}. The $\mathcal{L}_p$ norm is by definition
\begin{eqnarray}\label{19}&&\parallel \widetilde{\widetilde{K}}\left(X_{1},x,b\right)  \parallel_q=\left(\int\left(\widetilde{\widetilde{K}}\left(y,x,b\right) \right)^q f(y)dy\right)^{\frac{1}{q}}\nonumber\\\nonumber
&=&\left(\int\left(\prod\limits_{j=1}^{d}K_{{\rho(x_j,b_j)}}\left(t_{j}\right)\right)^q f(t_{1}^d)dt_{1}^{d}\right)^{\frac{1}{q}}
=\left(\mathsf{E}\left( f(\xi_{1}^d)\prod\limits_{j=1}^{d}K^{q-1}_{{\rho(x_j,b_j)}}\left(\xi_{j}\right)\right)\right)^{\frac{1}{q}},
\end{eqnarray}
where the product kernel $\prod\limits_{j=1}^{d}K_{{\rho(x_j,b_j)}}\left(\xi_{j}\right)$ is used as the pdf. $\{\xi_{j}\}$ are iid gamma distributed r.v.s with the  expectation $\mu_j=x_j$ and the variance $var = x_jb_j$. Analogously to the previous proofs, we use the Taylor expansion to find the expectation \eqref{19}, i.e. it holds
\begin{eqnarray}\label{20}&&\mathsf{E}\left(f(\xi_{1}^d)\prod\limits_{j=1}^{d}K^{q-1}_{{\rho(x_j,b_j)}}\left(\xi_{j}\right)\right)=
f(\mu_{1}^d)\prod\limits_{j=1}^{d}K^{q-1}_{{\rho(x_j,b_j)}}\left(\mu_{j}\right)\\\nonumber
&+&\sum\limits_{i=1}^{d}\frac{Var}{2}\frac{\partial^2}{\partial \xi_i^2}\left(f(\xi_{1}^d)\prod\limits_{j=1}^{d}K^{q-1}_{{\rho_1(x_j),b_j}}\left(\xi_{j}\right) \right)\Bigg|_{\xi=\mu}+o\left(\sum\limits_{j=1}^d b_j^2\right)\\\nonumber
&=&\left(\prod\limits_{j=1}^{d}K^{q-1}_{{\rho(x_j,b_j)}}\left(x_{j}\right)\right)\Bigg(f(x_{1}^d)\frac{(3-q)}{2}+\sum\limits_{i=1}^d
\frac{q(q-1)b_i}{2x_i}f(x_{1}^d)\\\nonumber
&-&(q-1)b_i\frac{\partial f(x_{1}^d)}{\partial x_i}+\frac{b_ix_i}{2}\frac{\partial^2 f(x_{1}^d)}{\partial x_i^2}\Bigg)+o\left(\sum\limits_{j=1}^d b_j^2\right).
\end{eqnarray}
Using Stirling's formula
\begin{eqnarray*}\label{stirling}\Gamma(z)=\sqrt{\frac{2\pi}{z}}\left(\frac{z}{e}\right)^z\left(1+O\left(\frac{1}{z}\right)\right),
\end{eqnarray*}
we can  rewrite the product kernel function as
\begin{eqnarray*}&&\prod\limits_{j=1}^{d}K_{{\rho(x_j,b_j)}}\left(x_{j}\right)=
\prod\limits_{j=1}^{d}\frac{x_j^{\frac{x_j}{b_j}-1}\exp(-\frac{x_j}{b_j})}{b_j^{\frac{x_j}{b_j}}\Gamma(\frac{x_j}{b_j})}
=\prod\limits_{j=1}^{d}\frac{x_j^{-\frac{1}{2}}b_j^{-\frac{1}{2}}}{\sqrt{2\pi}(1+O(b_j/x_j))}.
\end{eqnarray*}
Its upper bound is given by
\begin{eqnarray}\label{K1}&&\prod\limits_{j=1}^{d}K^{q-1}_{{\rho(x_j,b_j)}}\left(x_{j}\right)\leq
\prod\limits_{j=1}^{d}\frac{1}{(2\pi x_j b_j)^{\frac{q-1}{2}}}.
\end{eqnarray}
Hence, substituting \eqref{K1} into \eqref{20} we can rewrite \eqref{19} as
\begin{eqnarray*}&&\mathsf{E}\left(f(\xi_{1}^d)\prod\limits_{j=1}^{d}K^{q-1}_{{\rho(x_j,b_j)}}\left(\xi_{j}\right)\right)=
\left(\prod\limits_{j=1}^{d}(2\pi x_j b_j)^{\frac{1-q}{2}}\right)\Bigg(f(x_{1}^d)\frac{(3-q)}{2}\\
&+&\sum\limits_{i=1}^d\frac{q(q-1)b_i}{2x_i}f(x_{1}^d)-(q-1)b_i\frac{\partial f(x_{1}^d)}{\partial x_i}+\frac{b_ix_i}{2}\frac{\partial^2 f(x_{1}^d)}{\partial x_i^2}\Bigg)+o\left(\sum\limits_{j=1}^d b_j^2\right).
\end{eqnarray*}
When $b_1=b_2=\ldots=b_d=b$ holds, we get
\begin{eqnarray*}&&\parallel \widetilde{\widetilde{K}}\left(X_{1},x,b\right)  \parallel_q=
b^{(\tau+1)\frac{1-q}{2q}}\left(\prod\limits_{j=1}^{d}(2\pi x_j)^{\frac{1-q}{2q}}\right)\Bigg(f(x_{1}^d)\frac{(3-q)}{2}\\
&+&b\sum\limits_{i=1}^d\frac{q(q-1)}{2x_i}f(x_{1}^d)-(q-1)\frac{\partial f(x_{1}^d)}{\partial x_i}+\frac{x_i}{2}\frac{\partial^2 f(x_{1}^d)}{\partial x_i^2}\Bigg)^{1/q}\Bigg)+o\left(b^2\right).
\end{eqnarray*}
Let $p=q$ and, substituting the latter norm into Davydov's inequality, we get
 \begin{eqnarray*}&&\Big|Cov\left(\widetilde{\widetilde{K}}\left(X_{1},x,b\right),\widetilde{\widetilde{K}}\left(X_{1+k},x,b\right)\right)\Big|\leq 2\pi \alpha(k)^{\frac{1}{r}}\Bigg(b^{d\frac{1-q}{q}}\left(\prod\limits_{j=1}^{d}(2\pi x_j)^{\frac{1-q}{q}}\right)\\
&\cdot&\Bigg(f(x_{1}^d)\frac{(3-q)}{2}+b\sum\limits_{i=1}^d(q-1)\left(
\frac{q}{2x_i}f(x_{1}^d)
-\frac{\partial f(x_{1}^d)}{\partial x_i}\right)\\
&+&\frac{x_i}{2}\frac{\partial^2 f(x_{1}^d)}{\partial x_i^2}\Bigg)^{\frac{2}{q}}\Bigg)\Bigg)+o\left(b^2\right).
\end{eqnarray*}
Taking $p=q=2+\delta$, $r=\frac{2+\delta}{\delta}$ it can be deduced that the covariance of the multivariate pdf estimate is given by
\begin{eqnarray*}&&Cov(\widehat{f}(x_{1}^{d}))=\frac{2}{n}\sum\limits_{k=1}^{n}\left(1-\frac{k}{n}\right)Cov\left(\widetilde{\widetilde{K}}\left(X_{1},x,b\right),\widetilde{\widetilde{K}}_{1}\left(X_{1+k},x,b\right)\right)\\
&\leq&\frac{2b^{-d\frac{\delta+1}{\delta+2}}}{n}(2\pi)^{1-d\frac{\delta+1}{\delta+2}}\Bigg(\prod\limits_{j=1}^{d}x_j^{-\frac{\delta+1}{\delta+2}}\Bigg)\sum\limits_{k=1}^{n}\left(1-\frac{k}{n}\right)\alpha(k)^{\frac{\delta}{2+\delta}}\\
&\cdot&
\Bigg(f(x_{1}^d)\frac{1-\delta}{2}+b\sum\limits_{i=1}^d
\frac{(\delta+1)(\delta+2)}{2x_i}f(x_{1}^d)-(\delta+1)\frac{\partial f(x_{1}^d)}{\partial x_i}\\
&+&\frac{x_i}{2}\frac{\partial^2 f(x_{1}^d)}{\partial x_i^2}\Bigg)^{\frac{2}{q}}\Bigg)
+o\left(b^2\right).
\end{eqnarray*}
Let us introduce the following notations
\begin{eqnarray}\label{24}S(\delta,x_{1}^d)=\sum\limits_{i=1}^d
\frac{(\delta+1)(\delta+2)}{2x_i}f(x_{1}^d)-(\delta+1)\frac{\partial f(x_{1}^d)}{\partial x_i}+\frac{x_i}{2}\frac{\partial^2 f(x_{1}^d)}{\partial x_i^2},
\end{eqnarray}
\begin{eqnarray}\label{25}D(\delta,x_{1}^d)&=&2(2\pi)^{1-d\frac{\delta+1}{\delta+2}}\Bigg(\prod\limits_{j=1}^{d}x_j^{-\frac{\delta+1}{\delta+2}}\Bigg).
\end{eqnarray}
In these notations the covariance can be bounded by the following expression
\begin{eqnarray*}Cov(\widehat{f}(x_{1}^{d}))\leq\frac{D(\delta,x_{1}^d)}{n}b^{-d\frac{\delta+1}{\delta+2}}
\Bigg(f(x_{1}^d)\frac{1-\delta}{2}+bS(\delta,x_{1}^d)\Bigg)^{\frac{2}{2+\delta}}\!\!\int\limits_{1}^{\infty}\alpha(\tau)^{\frac{\delta}{2+\delta}}d\tau+o\left(b^2\right).
\end{eqnarray*}
Let us denote $\frac{\delta}{2+\delta}=\upsilon$, $0<\upsilon<1$. Then, we get the upper bound of the covariance \eqref{31}.
\subsection{Proof of Theorem \ref{thm2}}\label{Ap8}
Combining the bias, the variance and the covariance from Lemmas \ref{lem3} - \ref{lem7} we can write the MSE of the pdf estimate
\begin{eqnarray*}&&MSE(\widehat{f}(x_{1}^{d}))= Bias(\widehat{f}(x_{1}^{d}))^2+Var(\widehat{f}(x_{1}^{d}))+Cov(\widehat{f}(x_{1}^{d}))\\
&=&\left(\frac{b}{2}\sum\limits_{j=1}^{d}x_j\frac{\partial^2 f(x_{1}^{d})}{\partial x_j^2}\right)^2
-\frac{1}{n}\left(f(x_{1}^{d})+\frac{b}{2}\sum\limits_{j=1}^{d}x_j\frac{\partial^2 f(x_{1}^{d})}{\partial x_j^2}\right)^2\\
&+&\frac{1}{n}\left(\prod\limits_{j=1}^{d}\frac{b^{-\frac{d}{2}}x_j^{-1/2}}{2\sqrt{\pi}}\right)\left(f(x_{1}^d)+bv_1(x_{1}^{d})+b^2v_2(x_{1}^{d})\right)\\
&+&\frac{D(\upsilon,x_{1}^{d})}{n}b^{-d\frac{\upsilon+1}{2}}
\Bigg(bS(\upsilon,x_{1}^{d})+f(x_{1}^d)\frac{3\upsilon-1}{2(\upsilon-1)}\Bigg)^{1-\upsilon}\int\limits_{1}^{\infty}\alpha(\tau)^{\upsilon}d\tau+o\left(b^2\right).
\end{eqnarray*}
Since the rate of convergence in $b$ of the covariance is larger than one of the variance when $b\rightarrow0$ as $n\rightarrow\infty$, i.e.
\begin{eqnarray*}&&b^{-\frac{d}{2}}<b^{-\frac{d(\upsilon+1)}{2}},
\end{eqnarray*}
we can neglect the variance. Taking the  integral from the MSE and minimizing the resulted MISE in $b$ we get the optimal bandwidth \eqref{11}.
\subsection{Proof of Lemma \ref{lem77}}\label{Ap9}
\par Let us partition the covariance of the pdf estimate into two sums
 \begin{eqnarray}\label{21}Cov(\widehat{f}(x_{1}^{d}))&\leq&\frac{2}{n}\sum\limits_{\tau=2}^{c(n)}Cov\left(\widetilde{\widetilde{K}}\left(X_{1},x,b\right),\widetilde{\widetilde{K}}\left(X_{\tau},x,b\right)\right)
 \end{eqnarray}
\begin{eqnarray}&+&\frac{2}{n}\sum\limits_{\tau=c(n)+1}^{\infty}Cov\left(\widetilde{\widetilde{K}}\left(X_{1},x,b\right),\widetilde{\widetilde{K}}\left(X_{\tau},x,b\right)\right)=I_1+I_2.
\end{eqnarray}
Using Lemma \ref{lem5} the second sum in \eqref{21} can be bounded by
\begin{eqnarray*}I_2\leq\frac{D(\delta,x_{1}^{d})}{n}b^{-d\frac{\delta+1}{\delta+2}}
\Bigg(f(x_{1}^d)\frac{1-\delta}{2}+bS(\delta,x_{1}^{d})\Bigg)^{\frac{2}{2+\delta}}\sum\limits_{\tau=c(n)}^{\infty}\alpha(\tau)^{\frac{\delta}{2+\delta}}+o\left(b^2\right),
\end{eqnarray*}
where we use the notations \eqref{24} and \eqref{25}. Moreover, we can rewrite it as
\begin{eqnarray*}I_2&\leq&\frac{D(\delta,x_{1}^{d})}{c(n)n}b^{-d\frac{\delta+1}{\delta+2}}
\Bigg(f(x_{1}^d)\frac{1-\delta}{2}+bS(\delta,x_{1}^{d})\Bigg)^{\frac{2}{2+\delta}}\sum\limits_{\tau=2}^{\infty}\tau\alpha(\tau)^{\frac{\delta}{2+\delta}}+o\left(b^2\right).
\end{eqnarray*}
Let us denote $\delta=\frac{4\kappa}{1-2\kappa}$, $0<\kappa<\frac{1}{2}$. In this notations, we get
\begin{eqnarray*}I_2\leq\frac{D(\kappa,x_{1}^{d})}{c(n)n}b^{-d\frac{2\kappa+1}{2}}
\Bigg(f(x_{1}^d)\frac{6\kappa-1}{2(2\kappa-1)}+bS(\kappa,x_{1}^{d})\Bigg)^{1-2\kappa}\sum\limits_{\tau=2}^{\infty}\tau\alpha(\tau)^{2\kappa}+o\left(b^2\right).
\end{eqnarray*}
The first sum in \eqref{21} is determined and bounded by
\begin{eqnarray*} I_1&=&\frac{2}{n}\sum\limits_{\tau=2}^{c(n)}\Big|Cov\left(\widetilde{\widetilde{K}}\left(X_{1},x,b\right),\widetilde{\widetilde{K}}\left(X_{\tau},x,b\right)\right)\Big|\\
&=&\frac{2}{n}\sum\limits_{\tau=2}^{c(n)}\Big|\mathsf{E}(\widetilde{\widetilde{K}}\left(X_{1},x,b\right)\cdot \widetilde{\widetilde{K}}\left(X_{\tau},x,b\right)))-\mathsf{E}(\widetilde{\widetilde{K}}\left(X_{1},x,b\right))\cdot \mathsf{E}(\widetilde{\widetilde{K}}\left(X_{\tau},x,b\right))\Big|\\
&=&\frac{2}{n}\sum\limits_{\tau=2}^{c(n)}\Big|\int \limits_{0}^{\infty}\int\limits_{0}^{\infty} \widetilde{\widetilde{K}}\left(u,x,b\right) \widetilde{\widetilde{K}}\left(v,x,b\right)f_\tau(u,v)du dv\\
&-&\int\limits_{0}^{\infty} \widetilde{\widetilde{K}}\left(u,x,b\right)f(u)du  \int\limits_{0}^{\infty} \widetilde{\widetilde{K}}\left(v,x,b\right)f(v)dv\Big|
\end{eqnarray*}
\begin{eqnarray*}&=&\frac{2}{n}\sum\limits_{\tau=2}^{c(n)}\Big|\int \limits_{0}^{\infty}\int\limits_{0}^{\infty} \widetilde{\widetilde{K}}\left(u,x,b\right) \widetilde{\widetilde{K}}\left(v,x,b\right)(f_\tau(u,v)-f(u) f(v))du dv\Big|\\
&\leq&\frac{2}{n}\sum\limits_{\tau=2}^{c(n)}\int \limits_{0}^{\infty}\int\limits_{0}^{\infty} \Big|\widetilde{\widetilde{K}}\left(u,x,b\right) \widetilde{\widetilde{K}}\left(v,x,b\right)\Big||f_\tau(u,v)-f(u) f(v))|du dv.
\end{eqnarray*}
Under the condition of the lemma $|f_\tau(x,y)-f(x)f(y)|\leq M$ we get
\begin{eqnarray*}&&I_1\leq \frac{2M}{n}\sum\limits_{\tau=2}^{c(n)}\left[\int\limits_{0}^{\infty}\Big| \prod\limits_{j=1}^{d}K_{{\rho(x_j,b_j),b_j}}\left(u\right)\Big|du\right]^2=\frac{2Mc(n)}{n}.
\end{eqnarray*}
We aim to make the rate of convergence in $b$ of $I_1$ larger than the rate of $I_2$. Then it must hold
\begin{eqnarray*}c(n)< b^{\frac{-d(2\kappa+1)}{4}}.
\end{eqnarray*}
Let us choose, for example, $c(n)=b^{-\frac{d}{8}}$, $\kappa=1/4$. Then the rates of convergence of $I_1$ and $I_2$ are given by \eqref{32}.
\subsection{Proof of Theorem \ref{thm6}}\label{Ap10}
\par Using the bias, the variance and the covariance from Lemmas \ref{lem3}, \ref{lem4} and \ref{lem77}, let us write the MSE of the pdf estimate $\widehat{f}(x_{1}^{d})$ and find the optimal bandwidth $b$ that minimizes the MISE. The MSE is given by
\begin{eqnarray*}&&MSE(\widehat{f}(x_{1}^{d}))=\left(\frac{b}{2}\sum\limits_{j=1}^{d}x_j\frac{\partial^2 f(x_{1}^{d})}{\partial x_j^2}\right)^2
+\frac{1}{n}\left(b^{-\frac{d}{2}}\prod\limits_{j=1}^{d}\frac{x_j^{-1/2}}{2\sqrt{\pi}}\right)\\
&\cdot&\left(f(x_{1}^d)+bv_1(x_{1}^{d})+b^2v_2(x_{1}^{d})\right)-\frac{1}{n}\left(f(x_{1}^{d})+\frac{b}{2}\sum\limits_{j=1}^{d}x_j\frac{\partial^2 f(x_{1}^{d})}{\partial x_j^2}\right)^2\\
&+&\frac{2M}{nb^{\frac{d}{8}}}+\frac{D(\kappa,x_{1}^{d})}{nb^{\frac{d}{16}}}
\Bigg(f(x_{1}^d)\frac{6\kappa-1}{2(2\kappa-1)}+bS(\kappa,x_{1}^{d})\Bigg)^{1-2\kappa}\sum\limits_{\tau=2}^{\infty}\tau\alpha(\tau)^{2\kappa}+o\left(b^2\right).
\end{eqnarray*}
Since the rate of convergence in $b$ of the variance is larger than the rate of the covariance we can neglect the latter. Taking the integral from the MSE and minimizing it in $b$ we find the optimal bandwidth \eqref{33}.
\subsection{Proof of Lemma \ref{lem1}}\label{Ap1}
To find the bias of the estimate $\widehat{f}'_{x_k}(x_{1}^{d})$ we use the same technique as in \eqref{Ap3}.
The expectation is given by
\begin{eqnarray}\label{2}&&\mathsf{E}_{X}(\widehat{f}'_{x_k}(x_{1}^{d}))=\frac{1}{b_k}\mathsf{E}\left(\prod\limits_{j=1}^{d}K_{\rho(x_{j},b_j)}\left(X_{j}\right)L(X_{k},x_k,b_k)\right)\\
&=&\frac{1}{b_k}\int L(t_{k},x_k,b_k)f(t_{1}^d)\prod\limits_{j=1}^{d}K_{\rho(x_{j},b_j)}\left(t_{j}\right)dt_{1}^{d}
=\frac{1}{b_k}\mathsf{E}_{\xi}(f(\xi_{1}^d)L(\xi_{k},x_k,b_k))\nonumber\\\nonumber
&=&\frac{1}{b_k}\mathsf{E}_{\xi}(f(\xi_{1}^d)\ln(\xi_n))
-\frac{1}{b_k}\mathsf{E}_{\xi}(f(\xi_{1}^d)(\ln b_k + \Psi(\rho(x_k,b_k)))+o\left(\sum\limits_{j=1}^db_j\right),
\end{eqnarray}
where $\{\xi_j\}$ are iid  gamma distributed r.v.s with the expectation $\mu_{j}=x_j$ and the variance
$x_jb_j$.
Using the Taylor expansion in $\mu_{1}^{d}$, we can write
\begin{eqnarray*}&&f(\xi_{1}^d)\ln(\xi_n)=f(\mu_{1}^{d})\ln(\mu_n)
+\sum\limits_{j=1}^{d}(\xi_j-\mu_j)\frac{\partial \left(f(\xi_{1}^d)\ln(\xi_n)\right)}{\partial x_j}\Big|_{\mu}\\
&+&\frac{1}{2}\sum\limits_{j=1}^{d}(\xi_j-\mu_j)^2\frac{\partial^2 \left(f(\xi_{1}^d)\ln(\xi_n)\right)}{\partial x_j^2}\Big|_{\mu}\\
&+&\sum\limits_{j\neq l}(\xi_j-\mu_j)(\xi_l-\mu_l)\frac{\partial^2 \left(f(\xi_{1}^d)\ln(\xi_n)\right))}{\partial x_j\partial x_l}\Big|_{\mu}+o\left(\sum\limits_{j=1}^db_j\right).
\end{eqnarray*}
Taking the expectation from the right- and the left-hand sides of the latter we get the first term in \eqref{2}
\begin{eqnarray*}\mathsf{E}_{\xi}(f(\xi_{1}^d)\ln(\xi_n))&=&f(x_{1}^{d})\ln(x_k)
+\frac{1}{2}\sum\limits_{j=1}^{d}x_jb_j\frac{\partial^2 f}{\partial x_j^2}\ln(x_k)\\
&+&\frac{1}{2}x_kb_k\left(\frac{2}{x_k}\frac{\partial f(x_{1}^{d})}{\partial x_k}-\frac{f(x_{1}^{d})}{x^2_k}\right)+o\left(\sum\limits_{j=1}^db_j\right).
\end{eqnarray*}
For the second term in \eqref{2} we use  the approximation of the Digamma function when $\rho\rightarrow\infty$
\begin{eqnarray*}\Psi(\rho) &=& \ln \rho - \frac{1}{2\rho}-\frac{1}{12\rho^2}+\frac{1}{120\rho^4}-
\frac{1}{252\rho^6}+O\left(\frac{1}{\rho^8}\right).
\end{eqnarray*}
Hence, we can write that
\begin{eqnarray*}\ln b_k + \Psi(\rho(x_k,b_k)) &=&\ln x_k - \frac{b_k}{2x_k}-\frac{b_k^2}{12x_k^2}+o(b^2_n).
\end{eqnarray*}
Hence \eqref{2} can be rewritten as follows
\begin{eqnarray}\label{3}&&\mathsf{E}_{X}(\widehat{f}'_{x_k}(x_{1}^{d}))=f(x_{1}^{d})\frac{b_k}{12x_k^2}+\frac{\partial f(x_{1}^{d})}{\partial x_k}\\\nonumber
&+&\frac{1}{2}\sum\limits_{j=1}^{d}x_jb_j\frac{\partial^2 f(x_{1}^{d})}{\partial x_j^2}\left(\frac{1}{2x_k}+\frac{b_k}{12x_k^2}\right)+o\left(\sum\limits_{j=1}^db_j\right).
\end{eqnarray}
The multivariate expectation \eqref{3} coincides with its univariate version in \citet{DobrovidovMarkovich:13a}.
From \eqref{3} we get the bias of the pdf derivative estimate
\begin{eqnarray*}
Bias(\widehat{f}'_{x_k}(x_{1}^{d}))=f(x_{1}^{d})\frac{b_k}{12x_k^2}+\frac{1}{2}\sum\limits_{j=1}^{d}x_jb_j\frac{\partial^2 f(x_{1}^{d})}{\partial x_j^2}\left(\frac{1}{2x_k}+\frac{b_k}{12x_k^2}\right)+\!o\left(\sum\limits_{j=1}^db_j\right).
\end{eqnarray*}
When $b_1=b_2=\ldots=b_d=b$ holds, it follows that the bias of the pdf derivative is given by \eqref{34}.
\subsection{Proof of Lemma \ref{lem2}}\label{Ap2}
By definition, the variance for the estimate $\widehat{f}'_{x_k}(x_{1}^{d})$ is determined as
\begin{eqnarray}\label{9}
&& Var(\widehat{f}'_{x_k}(x_{1}^{d}))=\frac{1}{n}Var\left(\frac{1}{b_k}L(X_{k},x_k,b_k)\prod\limits_{j=1}^{d}K_{\rho(x_{j},b_j)}\left(X_{j}\right)\right)\\\nonumber
&=&\frac{1}{nb_k^2}\mathsf{E}\left(L^2(X_{k},x_k,b_k)\prod\limits_{j=1}^{d}K^2_{\rho(x_{j},b_j)}\left(X_{j}\right)\right)\\\nonumber
&-&\frac{1}{nb_k^2}\mathsf{E}^2\left(L(X_{k},x_k,b_k)\prod\limits_{j=1}^{d}K_{\rho(x_{j},b_j)}\left(X_{j}\right)\right).
\end{eqnarray}
The second term of the right-hand side of \eqref{9} is the square degree of \eqref{2}.  Hence, using \eqref{3} we can immediately write  that
\begin{eqnarray}\label{4}&&\frac{1}{n}\mathsf{E}_{X}^2(\widehat{f}'_{x_k}(x_{1}^{d}))
=\frac{1}{n}\Bigg(Bias^2(\widehat{f}'_{x_k}(x_{1}^{d}))+\left(\frac{\partial f(x_{1}^{d})}{\partial x_k}\right)^2\\\nonumber
&+&2\frac{\partial f(x_{1}^{d})}{\partial x_k}Bias(\widehat{f}'_{x_k}(x_{1}^{d}))+o\left(\sum\limits_{j=1}^db_j^2\right)\Bigg)\nonumber
\end{eqnarray}
hold. The first term in \eqref{9} can be represented as
\begin{eqnarray}\label{26}&&\frac{1}{nb_k^2}\mathsf{E}\left(L^2(X_{k},x_k,b_k)\prod\limits_{j=1}^{d}K^2_{\rho(x_{j},b_j)}\left(X_{j}\right)\right)\\\nonumber
&=&\frac{1}{nb_k^2}\int L^2(t_{k},x_k,b_k)f(t_{1}^d)\prod\limits_{j=1}^{d}K^2_{\rho(x_{j},b_j)}\left(t_{j}\right)dt_{1}^{d}\nonumber\\
&=&\frac{1}{nb_k^2}\int\left(\prod\limits_{j=1}^{d}\frac{t_j^{\frac{2x_j}{b_j}-2}\exp\left(\frac{-2t_j}{b_j}\right)}{b_j^{\frac{2x_j}{b_j}}\Gamma^2\left(\frac{x_j}{b_j}\right)}\right)
L^2(t_{k},x_k,b_k)f(t_{1}^d)dt_{1}^{d}.
\end{eqnarray}
Using the property of the gamma function
$\Gamma^2(\frac{x}{b}+1)=(\frac{x}{b})^2\Gamma^2(\frac{x}{b})$, \eqref{26} can be rewritten as
\begin{eqnarray}\label{10}&&\frac{1}{nb_k^2}\int\left(\prod\limits_{j=1}^{d}\frac{b_j^{-3}x_j^2\Gamma(\frac{2x_j}{b_j}-1)}{2^{\frac{2x_j}{b_j}-1}
\Gamma^2(\frac{x_j}{b_j}+1)}\frac{t_j^{\frac{2x_j}{b_j}-2}\exp(\frac{-2t_j}{b_j})}{\left(\frac{b_j}{2}\right)^{\frac{2x_j}{b_j}-1}\Gamma^(\frac{2x_j}{b_j}-1)}\right)
L^2(t_{k},x_k,b_k)f(t_{1}^d)dt_{1}^{d}\nonumber\\\nonumber
&=&\frac{1}{nb_k^2}\left(\prod\limits_{j=1}^{d}B(x_j,b_j)\right)\mathsf{E}_{\eta}\left(L^2(\eta_{k},x_k,b_k)f(\eta_{1}^d)\right)\\\nonumber
&=&\frac{1}{nb_k^2}\left(\prod\limits_{j=1}^{d}B(x_j,b_j)\right)\Bigg(\mathsf{E}_{\eta}\left(f(\eta_{1}^d)\ln^2\eta_k\right)
+\left(\ln b_k + \Psi\left(\frac{x_k}{b_k}\right)\right)^2\mathsf{E}_{\eta}\left(f(\eta_{1}^d) \right)\\
&-&2\left(\ln b_k + \Psi\left(\frac{x_k}{b_k}\right)\right)\mathsf{E}_{\eta}\left(f(\eta_{1}^d)\ln\eta_k\right)
\Bigg),
\end{eqnarray}
where $\{\eta_j\}$ are iid  gamma distributed r.v.s with the expectation $\mu_{j}=x_j-\frac{b_j}{2}$ and the variance
$\frac{x_jb_j}{2}-\frac{b_j^2}{4}$, and $B(x_j,b_j)$ is determined by \eqref{17}. To determine the expectations in \eqref{10} we use the same technique as  before.
Namely, the expectations are given by 
\begin{eqnarray*}\mathsf{E}_{\eta}(f(\eta_{1}^d)\ln(\eta_k))
&=&f(\mu_{1}^d)\left(\ln\mu_k-\frac{b_k}{4\mu_k}\right)+\sum\limits_{j=1}^{d}\left(\frac{x_jb_j}{2}-\frac{b_j^2}{4}\right)\frac{\partial^2 f(\mu_{1}^d)}{\partial x_j^2}\frac{\ln\mu_k}{2}\\
&+&\frac{b_k}{2}\frac{\partial f(\mu_{1}^d)}{\partial x_k}+o\left(\sum\limits_{j=1}^db_j^2\right),
\end{eqnarray*}
\begin{eqnarray*}&&\mathsf{E}_{\eta}(f(\eta_{1}^d))=f(\mu_{1}^d)+\frac{1}{2}\sum\limits_{j=1}^{d}\left(\frac{x_jb_j}{2}-\frac{b_j^2}{4}\right)\frac{\partial^2 f(\mu_{1}^d)}{\partial x_j^2}+o\left(\sum\limits_{j=1}^db_j^2\right),
\end{eqnarray*}
\begin{eqnarray*}\mathsf{E}_{\eta}(f(\eta_{1}^d)\ln^2(\eta_k))
&=&f(\mu_{1}^d)\left(\ln^2(\mu_k)+\frac{b_k}{2\mu_k}\left(1-\ln\mu_k\right)\right)+b_k\ln\mu_k\frac{\partial f(\mu_{1}^d)}{\partial x_k}\\
&+&\frac{\ln^2\mu_k}{2}\sum\limits_{j=1}^{d}\left(\frac{x_jb_j}{2}-\frac{b_j^2}{4}\right)\frac{\partial^2 f(\mu_{1}^d)}{\partial x_j^2}+o\left(\sum\limits_{j=1}^db_j^2\right).
\end{eqnarray*}
Using the known Tailor expansions
\begin{eqnarray*}\ln \mu_k&=& \ln x_k -\frac{b_k}{2x_k}-\frac{b_k^2}{8x_k^2}+o(b_k^2),\\
\frac{1}{\mu_k}&=&\frac{1}{x_k}+\frac{b_k}{2x_k^2}+\frac{b_k^2}{4x_k^3}+o(b_k^2)
\end{eqnarray*}
we can rewrite \eqref{10} as follows
\begin{eqnarray*}&&\frac{1}{n}\left(\prod\limits_{j=1}^{d}B(x_j,b_j)\right)
\Bigg(f(\mu_{1}^d)\Bigg(\frac{1}{4x^2_k}+\frac{1}{2x_kb_k}+\frac{7b_k}{48x_k^3}+\frac{7b^2_n}{576x_k^4}+\frac{b_k^3}{192x_k^5}\Bigg)\\
&+&\sum\limits_{j=1}^{d}\left(\frac{x_jb_j}{2}-\frac{b_j^2}{4}\right)\frac{\partial^2 f(\mu_{1}^d)}{\partial x_j^2}
\frac{b_k^2}{1152x_k^4}-\frac{b_k}{24x_k^2}\frac{\partial f(\mu_{1}^d)}{\partial x_k}+o\left(\sum\limits_{j=1}^db_j^2\right)
\Bigg).
\end{eqnarray*}
Therefore, using a Taylor expansion of the arguments of the functions $f(\mu_{1}^d)$ into $x_{1}^d$,
substituting it in latter and using \eqref{4}, we can write the variance of the pdf derivative estimate as
\begin{eqnarray*}&&Var(\widehat{f}'_{x_k}(x_{1}^{d}))=\\
&=&\frac{1}{n}\left(\prod\limits_{j=1}^{d}\frac{b_j^{-1/2}x_j^{-1/2}}{2\sqrt{\pi}}\right)
\Bigg(-\frac{b_k}{24x_k^2}\left(\frac{\partial f(x_{1}^d)}{\partial x_k}-\frac{b_j}{2}\frac{\partial^2 f(x_{1}^d)}{\partial x_k\partial x_j}+\sum\limits_{j=1}^{d}\frac{b_j^2}{4}\frac{\partial^3 f(x_{1}^d)}{\partial x_k\partial x_j^2}\right)\\\nonumber
&+&\frac{b_k^2}{1152x_k^4}\sum\limits_{j,i=n-\tau}^{n}\left(\frac{x_jb_j}{2}-\frac{b_j^2}{4}\right)
\left(\frac{b_i^2}{4}\frac{\partial^4 f(x_{1}^d)}{\partial x_j^2\partial x_i^2}-\frac{b_i}{2}\frac{\partial^3 f(x_{1}^d)}{\partial x_j^2\partial x_i}\right)+\frac{\partial^2 f(x_{1}^d)}{\partial x_j^2}\\\nonumber
&\cdot&\Bigg(\left(\frac{x_jb_j}{2}-\frac{b_j^2}{4}\right)\frac{b_k^2}{1152x_k^4}
+\frac{b_j^2}{4}\left(\frac{1}{4x^2_k}+\frac{1}{2x_kb_k}+\frac{7b_k}{48x_k^3}+\frac{7b^2_n}{576x_k^4}+\frac{b_k^3}{192x_k^5}\right)\Bigg)\\\nonumber
&+&\left(f(x_{1}^d)-\frac{b_j}{2}\frac{\partial f(x_{1}^d)}{\partial x_j}\right)
\Bigg(\frac{1}{4x^2_k}+\frac{1}{2x_kb_k}+\frac{7b_k}{48x_k^3}+\frac{7b^2_n}{576x_k^4}+\frac{b_k^3}{192x_k^5}\Bigg)\\\nonumber
&-&\frac{1}{n}\left(Bias^2(\widehat{f}'_{x_k}(x_{1}^{d}))+\left(\frac{\partial f(x_{1}^{d})}{\partial x_k}\right)^2+2\frac{\partial f(x_{1}^{d})}{\partial x_k}Bias(\widehat{f}'_{x_k}(x_{1}^{d}))\right)+o\left(\sum\limits_{j=1}^db_j^2\right).
\end{eqnarray*}
If $b_1=b_2=\ldots=b_d=b$ holds,  we get the variance of the pdf derivative estimate \eqref{38}.
\subsection{Proof of Lemma \ref{lem5}}\label{Ap5}
Now we apply Davydov's inequality for the covariance of the pdf derivative estimate
\begin{eqnarray*}&&|Cov\left(\widetilde{K}\left(X_{1},x,b\right),\widetilde{K}\left(X_{1+i},x,b\right)\right)|\leq\\
 &&2\pi \alpha(i)^{1/r}\parallel \widetilde{K}\left(X_{1},x,b\right) \parallel_q \parallel \widetilde{K}\left(X_{1+i},x,b\right) \parallel_p,
\end{eqnarray*}
where $p^{-1}+q^{-1}+r^{-1}=1$, $1\leq p, q, r\leq\infty$, $\alpha(i)$ is the mixing coefficient \eqref{1}. The $\mathcal{L}_p$ norm is given by
\begin{eqnarray}\label{27}&&\parallel \widetilde{K}\left(X_{1},x,b\right)  \parallel_q=\left(\int\left(\widetilde{K}\left(y,x,b\right) \right)^q f(y)dy\right)^{1/q}\nonumber\\
&=&\left(\int\left(L(t_{k},x_k,b_k)\prod\limits_{j=1}^{d}K_{{\rho(x_j,b_j),b_j}}\left(t_{j}\right)\right)^q f(t_{1}^d)dt_{1}^{d}\right)^{1/q}\nonumber\\
&=&\left(\int\left(L^q(t_{k},x_k,b_k) f(t_{1}^d)\prod\limits_{j=1}^{d}K^{q-1}_{{\rho(x_j,b_j),b_j}}\left(t_{j}\right) \right)\prod\limits_{j=1}^{d}K_{{\rho(x_j,b_j),b_j}}\left(t_{j}\right)dt_{1}^{d}\right)^{1/q}\nonumber\\
&=&\left(\mathsf{E}\left(L(\xi_{k},x_k,b_k)^q f(\xi_{1}^d)\prod\limits_{j=1}^{d}K^{q-1}_{{\rho(x_j,b_j),b_j}}\left(\xi_{j}\right)\right)\right)^{1/q},
\end{eqnarray}
where the kernel $\prod\limits_{j=1}^{d}K_{{\rho_1(x_j),b_j}}\left(\xi_{j}\right)$ is used as the pdf and $\{\xi_{j}\}$ are  iid gamma distributed  r.v.s with the expectation  $\mu_j=x_j$ and the variance $x_jb_j$. The Taylor expansion of \eqref{27} is the following
\begin{eqnarray*}&&L(\xi_{n},x_k,b_k)^q f(\xi_{1}^d)\prod\limits_{j=1}^{d}K^{q-1}_{{\rho(x_j,b_j),b_j}}\left(\xi_{j}\right)=\\
&=&\sum\limits_{i=1}^{d}\left(\xi_i-\mu_i\right)\frac{\partial}{\partial\xi_i}\left(L(\xi_{k},x_k,b_k)^q f(\xi_{1}^d)\prod\limits_{j=1}^{d}K^{q-1}_{{\rho(x_j,b_j),b_j}}\left(\xi_{j}\right)\right)\Bigg|_{\xi=\mu}\\
&+&\sum\limits_{i\neq l}\frac{(\xi_i-\mu_i)(\xi_l-\mu_l)}{2}\frac{\partial^2}{\partial\xi_i^2}\left(L(\xi_{k},x_k,b_k)^q f(\xi_{1}^d)\prod\limits_{j=1}^{d}K^{q-1}_{{\rho(x_j,b_j),b_j}}\left(\xi_{j}\right)\right)\Bigg|_{\xi=\mu}\\
&+&\sum\limits_{i=1}^{d}\frac{\left(\xi_i-\mu_i\right)^2}{2}\frac{\partial^2}{\partial\xi_i^2}\left(L(\xi_{k},x_k,b_k)^q f(\xi_{1}^d)\prod\limits_{j=1}^{d}K^{q-1}_{{\rho(x_j,b_j),b_j}}\left(\xi_{j}\right)\right)\Bigg|_{\xi=\mu}\\
&+&L(\mu_{k},x_k,b_k)^q f(\mu_{1}^d)\prod\limits_{j=1}^{d}K^{q-1}_{{\rho(x_j,b_j),b_j}}\left(\mu_{j}\right)+o\left(\sum\limits_{j=1}^db_j^2\right).
\end{eqnarray*}
Taking the expectation from the latter expression and
using the known Taylor expansions
\begin{eqnarray*}\label{L1}L(x_k,x_k, b_k)&=&\ln x_k-\ln b_k-\Psi\left(\frac{x_k}{b_k}\right)=
\frac{b_k}{2x_k}+\frac{b_k^2}{12x_k^2}+o(b_k^2),
\end{eqnarray*}
\begin{eqnarray*}L^q(x_k,x_k, b_k)&=&\left(\frac{b_k}{2x_k}\right)^q\left(1+q\frac{b_k}{6x_k}+\frac{q(q-1)}{2}\left(\frac{b_k}{6x_k}\right)^2\right),
\end{eqnarray*}
\begin{eqnarray*}\frac{\partial}{\partial \xi_k}L(\xi_n,x_k, b_k)&=&\frac{1}{\xi_k}\bigg|_{\xi_k=x_k}=\frac{1}{x_k},\quad
\frac{\partial^2}{\partial \xi_k^2}L(\xi_k,x_k, b_k)=-\frac{1}{\xi_n}\bigg|_{\xi_k=x_k}=-\frac{1}{x_k^2},
\end{eqnarray*}
the expectation \eqref{27} can be rewritten as
\begin{eqnarray*}&&\mathsf{E}\left(\prod\limits_{j=1}^{d}K^{q-1}_{{\rho_1(x_j),b_j}}\left(\xi_{j}\right)L_1(\xi_{n},x_k,b_k)^q f(\xi_{1}^d)\right)\\
&=&\left(\prod\limits_{j=1}^{d}K^{q-1}_{{\rho_1(x_j),b_j}}\left(x_{j}\right)\right)\Bigg(\left(\frac{b_k}{2x_k}\right)^q\left(1+q\frac{b_k}{6x_k}+\frac{q(q-1)}{2}\left(\frac{b_k}{6x_k}\right)^2\right) f(x_{1}^d)\\
&+&
\left(\frac{b_k}{2x_k}\right)^{q-1}\sum\limits_{i=1}^d\frac{x_ib_i}{2}\Bigg(\left(\frac{b_k}{2x_k}+\frac{qb_k^2}{12x_k^2}+\frac{q(q-1)b_k^3}{144x_k^3}\right)\\
&\cdot&\Bigg(f(x_{1}^d)\left(\frac{1-q}{b_ix_i}-\frac{q(1-q)}{x_i^2}\right)+\frac{\partial f(x_{1}^d)}{\partial x_i}\frac{2(1-q)}{x_i}
+\frac{\partial^2 f(x_{1}^d)}{\partial x_i^2}\Bigg)\\
&-&\left(1+\frac{(q-1)b_k}{6x_k}+\frac{(q-1)(q-2)b_k^2}{72x_k^2}\right)\left(f(x_{1}^d)\frac{q^2}{x_k^2}-\frac{2q}{x_k}\frac{\partial f(x_{1}^d)}{\partial x_k}\right)\Bigg)\\
&+&o\left(\sum\limits_{j=1}^db_j^2\right).
\end{eqnarray*}
If $b_1=b_2=\ldots=b_d=b$ holds and, using \eqref{K1}, we can obtain
 \begin{eqnarray*}&&\!\!\parallel \widetilde{K}\left(X_{1},x,b\right)  \parallel_q
\leq\left(\prod\limits_{j=1}^{d}(2\pi x_j)^{\frac{1-q}{2q}}\right)\frac{b^{d\frac{1-q}{2q}+1}}{2x_k}
\Bigg(\sum\limits_{i=1}^d
\left(1+\frac{qb}{6x_k}+\frac{q(q-1)b^2}{72x_k^2}\right)\\
&\cdot&\!\!\left(f(x_{1}^d)\frac{3-q}{2}+b(1-q)\frac{\partial f(x_{1}^d)}{\partial x_i}+
\frac{x_ib}{2}\frac{\partial^2 f(x_{1}^d)}{\partial x_i^2}-\frac{q(1-q)b}{2x_i}f(x_{1}^d)\right)\\
&+&\!\!x_i\left(1+\frac{(q-1)b}{6x_k}+\frac{(q-1)(q-2)b^2}{72x_k^2}\right)\left(2q\frac{\partial f(x_{1}^d)}{\partial x_k}-\frac{q^2}{x_k}f(x_{1}^d)\right)\Bigg)^{\frac{1}{q}}\!+\!o\left(b^2\right).
\end{eqnarray*}
The same is valid for $\parallel \widetilde{K}\left(X_{1+i},x,b\right) \parallel_p$. Let us take $p=q$. Hence, the covariance of the pdf derivative estimate is determined by
 \begin{eqnarray*}&&|Cov\left(\widetilde{K}\left(X_{1},x,b\right),\widetilde{K}\left(X_{1+k},x,b\right)\right)|\\
 &\leq& 2\pi \alpha(k)^{1/r}\Bigg(\left(\prod\limits_{j=1}^{d}(2\pi x_j)^{\frac{1-q}{q}}\right)\frac{b^{d\frac{1-q}{q}+2}}{4x_k^2}
 \Bigg(\sum\limits_{i=1}^d
\left(1+\frac{qb}{6x_k}+\frac{q(q-1)b^2}{72x_k^2}\right)\\
&\cdot&\left(f(x_{1}^d)\frac{3-q}{2}+b(1-q)\frac{\partial f(x_{1}^d)}{\partial x_i}+
\frac{x_ib}{2}\frac{\partial^2 f(x_{1}^d)}{\partial x_i^2}-\frac{q(1-q)b}{2x_i}f(x_{1}^d)\right)\\
&+&x_i\left(1+\frac{(q-1)b}{6x_k}+\frac{(q-1)(q-2)b^2}{72x_k^2}\right)\left(2q\frac{\partial f(x_{1}^d)}{\partial x_k}-\frac{q^2}{x_k}f(x_{1}^d)\right)\Bigg)^{\frac{2}{q}}\!+\!o\left(b^2\right).
\end{eqnarray*}
Taking $p=q=2+\delta$, $r=\frac{2+\delta}{\delta}$
and using the following notations
 \begin{eqnarray*}&&\!\!V(\delta,x_{1}^{d})=\sum\limits_{i=1}^d\bigg(-(\delta+1)\frac{\partial f(x_{1}^d)}{\partial x_i}+\frac{\delta(\delta+1)(\delta+2)}{36x_k^2}\frac{\partial f(x_{1}^d)}{\partial x_k}+
\frac{x_i}{2}\frac{\partial^2 f(x_{1}^d)}{\partial x_i^2}\\
 &-&\bigg(\frac{(\delta-1)(\delta+1)(\delta+2)}{144x_k^2}-\frac{(\delta+1)(\delta+2)}{2x_i}+\frac{\delta(\delta+1)(\delta+2)^2}{72x_k^3}\bigg)f(x_{1}^d)
\Bigg),
 \end{eqnarray*}
  \begin{eqnarray*}W(\delta,x_{1}^{d})&=&\sum\limits_{i=1}^d\bigg(\bigg(\frac{1-\delta}{4x_k}+\frac{(\delta+1)(\delta+2)}{2x_i}-\frac{x_i(\delta+1)(\delta+2)^2}{6x_k^2}\bigg)f(x_{1}^d)\\
&-&
(\delta+1)\frac{\partial f(x_{1}^d)}{\partial x_i}+\frac{(\delta+1)(\delta+2)x_i}{3x_k}\frac{\partial f(x_{1}^d)}{\partial x_k}+\frac{x_i}{2}\frac{\partial^2 f(x_{1}^d)}{\partial x_i^2}\bigg),
 \end{eqnarray*}
  \begin{eqnarray*}L(\delta,x_{1}^{d})&=&\sum\limits_{i=1}^d\Bigg(f(x_{1}^d)\frac{(1-\delta)}{2}+x_i\bigg(2(2+\delta)\frac{\partial f(x_{1}^d)}{\partial x_k}-\frac{(2+\delta)^2}{x_k}f(x_{1}^d)\bigg)\Bigg),
   \end{eqnarray*}
   \begin{eqnarray*}R(\delta,x_{1}^{d})&=&\left(\prod\limits_{j=1}^{d}x_j^{-\frac{\delta+1}{\delta+2}}\right)\frac{(2\pi)^{\frac{\delta+2-d(\delta+1)}{\delta+2}}}{2x_k^2},
   \end{eqnarray*}
the covariance of the pdf derivative estimate is given by
 \begin{eqnarray*}&&C(\widehat{f}'_{x_k}(x_{1}^{d})=\frac{2}{nb^2}\sum\limits_{i=1}^{n}\left(1-\frac{i}{n}\right)Cov\left(\widetilde{K}\left(X_{1},x,b\right),\widetilde{K}\left(X_{1+i},x,b\right)\right)\\
&\leq&\frac{R(\delta,x_{1}^{d})}{b^{d\frac{\delta+1}{\delta+2}}n}
\Bigg(b^2V(\delta,x_{1}^{d})+bW(\delta,x_{1}^{d})+L(\delta,x_{1}^{d})\Bigg)^{\frac{2}{\delta+2}}\int\limits_{1}^{\infty}\alpha(\tau)^{\frac{\delta}{2+\delta}}d\tau+o\left(b^2\right).
\end{eqnarray*}
Let us denote $\frac{\delta}{2+\delta}=\upsilon$, $0<\upsilon<1$. In these notations, we get \eqref{cov}.
\subsection{Proof of Theorem \ref{thm}}\label{Ap6}
Using the bias, the variance and the covariance in Lemmas \ref{lem1}, \ref{lem2}, \ref{lem5}, we can write the MSE for the estimate of the pdf derivative
\begin{eqnarray*}&&MSE(\widehat{f}'_{x_k}(x_{1}^{d}))= (Bias\widehat{f}'_{x_k}(x_{1}^{d}))^2+Var(\widehat{f}'_{x_k}(x_{1}^{d}))+Cov(\widehat{f}'_{x_k}(x_{1}^{d}))\\
&=&\Bigg(bB_1(x_{1}^{d})+b^2B_2(x_{1}^{d})\Bigg)^2\left(1-\frac{1}{n}\right)\\
&+&\frac{b^{-\frac{d}{2}}}{n}\left(\prod\limits_{j=1}^{d}\frac{x_j^{-1/2}}{2\sqrt{\pi}}\right)
\Bigg(bV_1(x_{1}^{d})+b^2V_2(x_{1}^{d})+\frac{1}{b}V_3(x_{1}^{d})+V_4(x_{1}^{d})\Bigg)\\
&-&\frac{1}{n}\left(\left(\frac{\partial f(x_{1}^{d})}{\partial x_k}\right)^2+2\frac{\partial f(x_{1}^{d})}{\partial x_k}\left(bB_1(x_{1}^{d})+b^2B_2(x_{1}^{d})\right)\right)
+\frac{R(\upsilon,x_{1}^{d})}{n}b^{-d\frac{\upsilon+1}{2}}\\
&\cdot&\Bigg(b^2V(\upsilon,x_{1}^{d})+bW(\upsilon,x_{1}^{d})+L(\upsilon,x_{1}^{d})\Bigg)^{1-\upsilon}\int\limits_{1}^{\infty}\alpha(\tau)^{\upsilon}d\tau+o\left(b^2\right).
\end{eqnarray*}
Since the rate of convergence of the covariance is less than that one of the variance, we can neglect the covariance term. Taking the integral from the MSE and minimizing the MISE in $b$ the optimal bandwidth \eqref{35} is obtained.
\begin{acknowledgements}
L.A. Markovich  was supported by the Russian Science Foundation grant (14-50-00150).
\end{acknowledgements}


\bibliographystyle{spbasic}      
\bibliography{reference}

\end{document}